\documentclass[draft]{article}
\usepackage{amsfonts}
\usepackage{graphicx}
\usepackage{amsmath}
\usepackage[english]{babel}
\usepackage{url}
\newcommand{\be}{\begin{equation}}
\newcommand{\ee}{\end{equation}}
\newcommand{\bt}{\begin{table}}
\newcommand{\et}{\end{table}}
\newcommand{\btr}{\begin{tabular}}
\newcommand{\etr}{\end{tabular}}
\newcommand{\bd}{\begin{displaymath}}
\newcommand{\ed}{\end{displaymath}}
\newcommand{\bc}{\begin{center}}
\newcommand{\ec}{\end{center}}
\newcommand{\bea}{\begin{eqnarray}}
\newcommand{\eea}{\end	{eqnarray}}
\newcommand{\ba}{\begin{array}}
\newcommand{\ea}{\end{array}}
\newcommand{\bi}{\begin{itemize}}
\newcommand{\ei}{\end{itemize}}
\newcommand{\ben}{\begin{enumerate}}
\newcommand{\een}{\end{enumerate}}
\newcommand{\bn}{\begin{note}}
\newcommand{\en}{\end{note}}
\newcommand{\bfg}{\begin{figure}}
\newcommand{\efg}{\end{figure}}
\newcommand{\bpm}{\begin{pmatrix}}
\newcommand{\epm}{\end{pmatrix}}
\newcommand{\bv}{\begin{verbatim}}
\newcommand{\ev}{\end{verbatim}}
\newcommand{\ds}{\displaystyle}
\newcommand{\bs}{\boldsymbol}
\newtheorem{theorem}{Theorem}
\newtheorem{assumption}[theorem]{Assumption}

\newtheorem{corollary}[theorem]{Corollary}

\newtheorem{definition}[theorem]{Definition}
\newtheorem{example}[theorem]{Example}
\newtheorem{exercise}[theorem]{Exercise}
\newtheorem{lemma}[theorem]{Lemma}

\newtheorem{problem}[theorem]{Problem}
\newtheorem{proposition}[theorem]{Proposition}
\newtheorem{remark}[theorem]{Remark}

\newtheorem{method}[theorem]{Method}
\newenvironment{proof}[1][Proof]{\textbf{#1.} }{\ \rule{0.5em}{0.5em}}
\newcommand{\bpb}{\begin{problem}}
\newcommand{\epb}{\end{problem}}
\newcommand{\bpf}{\begin{proof}}
\newcommand{\epf}{\end{proof}}
\newcommand{\bth}{\begin{theorem}}
\newcommand{\eth}{\end{theorem}}
\newcommand{\bpr}{\begin{proposition}}
\newcommand{\epr}{\end{proposition}}
\newcommand{\bl}{\begin{lemma}}
\newcommand{\el}{\end{lemma}}
\newcommand{\bdf}{\begin{definition}}
\newcommand{\edf}{\end{definition}}
\newcommand{\bco}{\begin{corollary}}
\newcommand{\eco}{\end{corollary}}
\newcommand{\br}{\begin{remark}}
\newcommand{\er}{\end{remark}}
\newcommand{\bex}{\begin{example}}
\newcommand{\eex}{\end{example}}
\newcommand{\bexr}{\begin{exercise}}
\newcommand{\eexr}{\end{exercise}}
\newcommand{\bme}{\begin{method}}
\newcommand{\eme}{\end{method}}
\newcommand{\bas}{\begin{assumption}}
\newcommand{\eas}{\end{assumption}}
\begin{document}
\title{Semigroup discretization and\\spectral approximation for\\linear nonautonomous delay differential equations}
\author{ Dimitri Breda$^{1}$, Stefano Maset$^{2}$ and Rossana Vermiglio$^{1}$\\
[.4em]{\small{\it $^{1}$Department of Mathematics and Computer Science}}\\
[-.3em]{\small{\it University of Udine}}\\
[-.3em]{\small{\it via delle Scienze 206, I33100 Udine - Italy}}\\
[-.3em]{\small{\it e-mail:} \url{{dimitri.breda},{rossana.vermiglio}@dimi.uniud.it}}\\
[.4em]{\small{\it $^{2}$Department of Mathematics and Computer Science}}\\
[-.3em]{\small{\it University of Trieste}}\\
[-.3em]{\small{\it via Valerio 12, I34127 Trieste - Italy}}\\
[-.3em]{\small{\it e-mail:} \url{maset@univ.trieste.it}}}
\date{\today}
\maketitle
\begin{abstract}
This paper deals with the approximation of the spectrum of linear and nonautonomous delay differential equations through the reduction of the relevant evolution semigroup from infinite to finite dimension. The focus is placed on classic collocation, even though the requirements that a numerical scheme has to fulfill in order to allow for a correct approximation of the spectral elements are recalled. This choice, motivated by the analyticity of the underlying eigenfunctions, allows for a convergence of infinite order, as rigorously demonstrated through {\it a priori} error bounds when Chebyshev nodes are adopted. Fundamental applications such as determination of asymptotic stability of equilibria (autonomous case) and limit cycles (periodic case) follow at once.
\end{abstract}
{\small {\it Keywords:} delay differential equations, spectrum, evolution semigroup, numerical collocation\\
MSC[2010]  34K30, 34L16, 65L03, 65L07, 65L15, 47D06}
\section{Introduction}\label{s_intro}
In the recent decades, Delay Differential Equations (DDEs), and more general functional differential equations, have attracted the attention of diverse scientific communities, beyond that of mathematicians, ranging from automatic control to physics, through population dynamics and bio-mathematics, to name a few. A central question from a dynamical point of view is that of being able to determine or foresee the asymptotic stability of equilibria and periodic orbits of nonlinear systems through suitable linearization. Despite the great effort, the utilization of well-established analytical results such as spectral bounds and Lyapunov stability turns out to be rather lacking from the practical point of view of applications and, at best, suitable for restricted sub-classes (e.g. single discrete delay, second order systems, etc.). A major difficulty is clearly manifested in the infinite dimensional nature of these models arising when their time evolution is described in a suitable state space.  As a natural consequence \cite[p.109]{hvl93}, a number of approximation techniques have been proposed, mostly based on computing the characteristic values (roots, multipliers, Lyapunov exponents) associated to the system, see e.g. \cite{bresisc05,bmbas04,engacm02,engsjna02,farmer82,ig02,jar08,vlr08,vz09}.

When investigating on the stability (but not only) of
\be\label{model}
x(t)=\mathcal{F}(t,x_{t}),
\ee
where $\mathcal{F}:[0,\infty)\times\mathcal{C}\rightarrow\mathbb{C}^{d}$ is linear with $\mathcal{C}$ the state space of elements $x_{t}\in\mathcal{C}$ defined as
\bd
x_{t}(\theta):=x(t+\theta),\;\theta\in[-\tau,0],
\ed
according to the standard Hale-Krasovskii notation \cite{hvl93,kra59}, $\tau$ being the maximum delay and $d$ the number of equations, the state space description of the model is advantageous, and the classic literature resorts to the Banach space of continuous functions $\mathcal{C}:=C(-\tau,0;\mathbb{C}^{d})$, \cite{belzen03,diekmann95,hvl93,wu96}. This choice seems to be motivated by the fact that, for rather general selections of the space of initial data, the ``smoothing effect'' \cite{belzen03} makes the (forward) solution be continuous anyway: ``{\it...if some other space than continuous functions is used for initial data, then the solution lies in $C$...Therefore, for the fundamental theory, the space of initial data does not play a role which is too significant.}'' \cite[p.33]{hvl93}. Anyway, Hale continues his comment by adding ``{\it However, in the applications, it is sometimes convenient to take initial data with fewer or more restrictions.}'' In this sense, an alternative which has been quite studied is represented by the Hilbert product space $\mathbb{C}\times L^{2}(-\tau,0;\mathbb{C}^{d})$, \cite{bddm92,bt69,dm72,hrs08,pei82}. This second choice is often justified in the context of quadratic feedback control and linear filtering for retarded systems \cite{del77,hrs08,vin78}, for approximation reasons \cite{ik91,kap86}, or when orthogonality is necessary \cite{brejmaa09}. In this manuscript the choice is that of continuous functions, reserving to present analogous arguments for the Hilbert space in forthcoming works of the authors, as already announced (and partially developed) in \cite{breifac210}.

Once the proper state space is chosen, the long-time behavior of the evolution can be determined through the knowledge of the spectrum of infinite dimensional maps such as the semigroup of solution operators and its generator in the autonomous case, the monodromy operator for periodic problems, the evolution family in the nonautonomous case. The reduction of such operators to finite dimension allows to consider standard eigenvalue problems which can be easily solved, hopefully providing accurate estimates for the stability indicators (e.g. the rightmost root, the dominant multiplier or the largest exponent). Usually, the construction of the finite dimensional approximation represents the ``easy'' step, although rather technical difficulties may arise when systems with possibly multiple discrete and distributed delays are the case. Instead, the theoretical (i.e. not experimental) analysis of convergence is far away to be as simple and direct as the implementation of the numerical scheme. Here both aspects are equally treated with rigor. In particular, as far as convergence is concerned, the general properties that a underlying numerical method has to possess in order to lead to a satisfactory approximation of the spectral elements (eigenvalues, mutliplicities, eigenspaces, etc.) are recalled. To this aim, the (different) theories on spectral approximation of linear operators developed in the monographs \cite{chat83,ggk90} are used as (alternative) background.

Beyond this general treatment, the focus is placed on the {\it pseudospectral} approach, namely {\it collocation} together with {\it polynomial interpolation}: the core consists in substituting the exact operation (e.g. differentiation) to be done on a given function over a selected distribution of nodes with the same operation as applied to the interpolating polynomial. The method, well-known for numerical integration of ordinary and partial differential equations as well as for the associated eigenvalue problems (see \cite{tref00} for a guide), in the context of DDEs was first presented in \cite{bresisc05} for the approximation of the spectrum of the infinitesimal generator for autonomous problems and in different guises also in \cite{breth04,breifac06,brenm09,bmbas04,gp06,vlr08}. It benefits from the infinite regularity of the underlying eigenfunctions of the system at hand \cite{hvl93} and, by choosing to operate on Chebyshev nodes, it is shown to be the unique one able to exploit all this regularity in problems with delay, too. As a result, the spectral elements are approximated at a convergence rate of infinite order, as it is rigorously demonstrated in the manuscript. To the best of the authors' knowledge, this represents the first (and only) complete proof of convergence in the field of DDEs.

The paper is structured as follows. The prototypical model for (\ref{model}) is introduced in Section \ref{s_model} together with its state space representation through the associated evolution family. After some preliminaries discussed in Section \ref{s_notC}, Section \ref{s_discrC} deals with the numerical discretization. The convergence analysis is tackled in Section \ref{s_convC}, precisely for the discretization scheme in Section \ref{s_conv1} and for the spectral elements in Sections \ref{s_conv2} and \ref{s_conv3}. In Section \ref{s_appl} applications are briefly discussed and commented. Appendix \ref{s_app} collects Definitions and Lemmas useful and necessary to prove the main results of Section \ref{s_convC}.
\section{Model, state space and evolution}\label{s_model}
As a prototypical model for (\ref{model}), we consider the scalar linear DDE with nonauto\-nomous coefficients
\be\label{dde}
x'(t)=a(t)x(t)+b(t)x(t-\tau)+\int\limits_{-\tau}^{0}c(t,\theta)x(t+\theta)d\theta
\ee
where $\tau>0$ is the maximum delay and $[0,+\infty)\ni t\mapsto a(t),b(t)\in\mathbb{C}$ and $[0,+\infty)\times[-\tau,0]\ni(t,\theta)\mapsto c(t,\theta)\in\mathbb{C}$ are continuous functions. All the arguments developed in the sequel apply as well to more general cases with matrix coefficients and multiple discrete or distributed delays, the extension concerning only technicalities useless to the treatment proposed in the paper.

In order to focus on a well-posed Initial Value Problem (IVP) for (\ref{dde}), it is necessary to specify a suitable set of initial data. The classic literature \cite{belzen03,diekmann95,hvl93,wu96} resorts to the Banach space of continuous functions $(\mathcal{C},\|\cdot\|_{\mathcal{C}})$ with
\bd
\mathcal{C}:=C(-\tau,0;\mathbb{C})
\ed
and
\be\label{normC}
\|\psi\|_{\mathcal{C}}:=\max\limits_{\theta\in[-\tau,0]}|\psi(\theta)|.
\ee
Then, given two reals $r\geq s$ and $\varphi\in\mathcal{C}$, the IVP
\be\label{IVPC}
\left\{
\setlength\arraycolsep{0.1em}\ba{ll}
\ds x'(t)=a(t)x(t)+b(t)x(t-\tau)+\int\limits_{-\tau}^{0}c(t,\theta)x(t+\theta)d\theta,&\;t\in[s,r]\\ 
\ds x(s+\theta)=\varphi(\theta),&\;\theta\in[-\tau,0]
\ea
\right.
\ee
admits the existence of a unique continuous solution $x\in[s-\tau,r]$, which continuously depends on the initial data \cite{hvl93}. This allows to introduce the linear and bounded operator $T(r,s):\mathcal{C}\rightarrow\mathcal{C}$ given by
\be\label{Tr}
T(r,s)x_{s}=x_{r},
\ee
$x_{r}$, the state of the system at time $r\geq s$.
In particular, $x_{s}=\varphi$. The two-parameters family $\{T(r,s)\}_{r\geq s}$ is a strongly continuous {\it evolution family} and it is eventually compact, i.e. $T(r,s)$ is compact as soon as $r\geq s+\tau$, see \cite{hvl93,diekmann95} and \cite{hrs08} for definition and properties.

In this work we are interested in the approximation of the {\it spectral elements} of $T(r,s)$, i.e. the nonzero eigenvalues as well as their multiplicities and eigenspaces, fundamental for stability purposes.
\section{Preliminaries and notation}\label{s_notC}
Depending on the role of the various mathematical objects, in general we use normal case for operators and functions (infinite dimension), bold case for matrices and vectors (finite dimension).

After the time translation $t\mapsto s+t$, set
\bd
\setlength\arraycolsep{0.1em}\ba{l}
\ds\mathcal{C}^{-}:=\mathcal{C}=C(-\tau,0;\mathbb{C})\\
\ds\mathcal{C}^{+}:=C(0,r_{s};\mathbb{C})\\
\ds\mathcal{C}^{\pm}:=C(-\tau,r_{s};\mathbb{C})
\ea
\ed
with $r_{s}:=r-s$ (the same convention holds for function spaces other than $\mathcal{C}$) and rewrite the IVP (\ref{IVPC}) as
\be\label{IVPG_{s}yC}
\left\{
\setlength\arraycolsep{0.1em}\ba{ll}
\ds y'(t)=(G_{s}y)(t),&\;t\in[0,r_{s}]\\
\ds y(\theta)=\varphi(\theta),&\;\theta\in[-\tau,0],
\ea
\right.
\ee
for $y(t):=x(s+t)$, $x$ solution of (\ref{IVPC}), where the operator $G_{s}:\mathcal{C}^{\pm}\rightarrow\mathcal{C}^{+}$ is defined as
\be\label{G_{s}}
(G_{s}y)(t):=a_{s}(t)y(t)+b_{s}(t)y(t-\tau)+\int\limits_{-\tau}^{0}c_{s}(t,\theta)y(t+\theta)d\theta,\;t\in[0,r_{s}],
\ee
where
\bd
a_{s}(t):=a(s+t),
\ed
\bd
b_{s}(t):=b(s+t)
\ed
and, for all $\theta\in[-\tau,0]$,
\bd
c_{s}(t,\theta):=c(s+t,\theta).
\ed
The solution $y\in\mathcal{C}^{\pm}$ in $[-\tau,r_{s}]$ is intended as divided into
\bd
\left\{
\setlength\arraycolsep{0.1em}\ba{ll}
\ds y^{+}(t):=y(t),&\;t\in[0,r_{s}]\\
\ds y^{-}(t):=y(t)=\varphi(t),&\;t\in[-\tau,0].
\ea
\right.
\ed

For a given positive integer $M$, consider the grid of distinct nodes
\bd
\Omega_{M}^{-}:=\{-\tau=:\theta_{M,M}^{-}<\cdots<\theta_{M,0}^{-}:=0\}
\ed
in $[-\tau,0]$ and set $\mathcal{C}_{M}^{-}:=\mathbb{C}^{M+1}$ as the discrete counterpart of $\mathcal{C}^{-}$, i.e. a function $f^{-}\in\mathcal{C}^{-}$ is discretized by the vector
\bd
\bs{f}_{M}^{-}=\mathcal{R}_{M}^{-}f^{-}=(f^{-}(\theta_{M,0}^{-}),\ldots,f^{-}(\theta_{M,M}^{-}))^{T}\in\mathcal{C}_{M}^{-}
\ed
where $\mathcal{R}_{M}^{-}:\mathcal{C}^{-}\rightarrow\mathcal{C}_{M}^{-}$ is the restriction operator associating to a function its grid values at the nodes $\Omega_{M}^{-}$. Let $\mathcal{P}_{M}^{-}:\mathcal{C}_{M}^{-}\rightarrow\Pi_{M}^{-}\subset\mathcal{C}^{-}$ be the prolongation operator defined as
\bd
(\mathcal{P}_{M}^{-}\bs{v}^{-})(t):=\sum\limits_{j=0}^{M}\ell_{M,j}^{-}(t)v^{-}_{j},\;t\in[-\tau,0],
\ed
for
\bd
\bs{v}^{-}=(v^{-}_{0},\ldots,v^{-}_{M})^{T}\in\mathcal{C}_{M}^{-}
\ed
and
\bd
\ell_{M,j}^{-}(t):=\prod_{\substack{k=0\\k\neq j}}^{M}\frac{t-\theta_{M,k}^{-}}{\theta_{M,j}^{-}-\theta_{M,k}^{-}},\;j=0,\ldots,M,
\ed
the Lagrange basis polynomials relevant to the nodes $\Omega_{M}^{-}$. Then $f_{M}^{-}=\mathcal{P}_{M}^{-}\bs{f}_{M}^{-}=\mathcal{P}_{M}^{-}\mathcal{R}_{M}^{-}f^{-}\in\Pi_{M}^{-}$ is the polynomial of degree at most $M$ interpolating $f^{-}$ at the nodes $\Omega_{M}^{-}$, in fact\bd
\mathcal{R}_{M}^{-}\mathcal{P}_{M}^{-}=\bs{I}_{M}^{-},
\ed
$\bs{I}_{M}^{-}:\mathcal{C}_{M}^{-}\rightarrow\mathcal{C}_{M}^{-}$ being the identity in $\mathcal{C}_{M}^{-}$,  while
\be\label{LM-}
\mathcal{P}_{M}^{-}\mathcal{R}_{M}^{-}=\mathcal{L}_{M}^{-},
\ee
$\mathcal{L}_{M}^{-}:\mathcal{C}^{-}\rightarrow\Pi_{M}^{-}\subset\mathcal{C}^{-}$ being the Lagrange interpolation operator on $\Omega_{M}^{-}$.

Similarly, for a given positive integer $N$, let
\be\label{grid+}
\Omega_{N}^{+}:=\{0<\theta_{N,1}^{+}<\cdots<\theta_{N,N}^{+}<r_{s}\}
\ee
be a grid of distinct nodes in $(0,r_{s})$, together with the auxiliary node $\theta_{N,0}^{+}:=0$ $(=\theta_{M,0}^{-})$ and set $\mathcal{C}_{N}^{+}:=\mathbb{C}^{N+1}$ as the discrete counterpart of $\mathcal{C}^{+}$, i.e. a function $f^{+}\in\mathcal{C}^{+}$ is discretized by the vector
\bd
\bs{f}_{N}^{+}=\mathcal{R}_{N,0}^{+}f^{+}=(f^{+}(\theta_{N,0}^{+}),\ldots,f^{+}(\theta_{N,N}^{+}))^{T}\in\mathcal{C}_{N}^{+}
\ed
where $\mathcal{R}_{N,0}^{+}:\mathcal{C}^{+}\rightarrow\mathcal{C}_{N}^{+}$ is the restriction operator associating to a function its grid values at the nodes $\{\theta_{N,0}^{+}\}\cup\Omega_{N}^{+}$. Let $\mathcal{P}_{N,0}^{+}:\mathcal{C}_{N}^{+}\rightarrow\Pi_{N}^{+}\subset\mathcal{C}^{+}$ be the prolongation operator defined as
\bd
(\mathcal{P}_{N,0}^{+}\bs{v}^{+})(t):=\sum\limits_{j=0}^{M}\ell_{N,j}^{+}(t)v^{+}_{j},\;t\in[0,r_{s}],
\ed
for
\bd
\bs{v}^{+}=(v^{+}_{0},\ldots,v^{+}_{N})^{T}\in\mathcal{C}_{N}^{+}
\ed
and
\bd
\ell_{N,j}^{+}(t):=\prod_{\substack{k=0\\k\neq j}}^{N}\frac{t-\theta_{N,k}^{+}}{\theta_{N,j}^{+}-\theta_{N,k}^{+}},\;j=0,\ldots,N,
\ed
the Lagrange basis polynomials relevant to the nodes $\{\theta_{N,0}^{+}\}\cup\Omega_{N}^{+}$. Then $f_{N}^{+}=\mathcal{P}_{N,0}^{+}\bs{f}_{N}^{+}=\mathcal{P}_{N,0}^{+}f^{+}\mathcal{R}_{N,0}^{+}\in\Pi_{N}^{+}$ is the polynomial of degree at most $N$ interpolating $f^{+}$ at the nodes $\{\theta_{N,0}^{+}\}\cup\Omega_{N}^{+}$, in fact
\bd
\mathcal{R}_{N,0}^{+}\mathcal{P}_{N,0}^{+}=\bs{I}_{N}^{+},
\ed
$\bs{I}_{N}^{+}:\mathcal{C}_{N}^{+}\rightarrow\mathcal{C}_{N}^{+}$ being the identity in $\mathcal{C}_{N}^{+}$, while
\bd
\mathcal{P}_{N,0}^{+}\mathcal{R}_{N,0}^{+}=\mathcal{L}_{N,0}^{+},
\ed
$\mathcal{L}_{N,0}^{+}:\mathcal{C}^{+}\rightarrow\Pi_{N}^{+}\subset\mathcal{C}^{+}$ being the Lagrange interpolation operator on $\{\theta_{N,0}^{+}\}\cup\Omega_{N}^{+}$. In the sequel it will be necessary to refer rather to the Lagrange interpolation operator on $\Omega_{N}^{+}$ (i.e. without $\theta_{N,0}^{+}$). This latter will be denoted by $\mathcal{L}_{N}^{+}:\mathcal{C}^{+}\rightarrow\Pi_{N-1}^{+}\subset\mathcal{C}^{+}$ and given as
\be\label{LN+}
\mathcal{L}_{N}^{+}=\mathcal{P}_{N}^{+}\mathcal{R}_{N}^{+}
\ee
with obvious meaning of the restriction and prolongation operators relevant to $\Omega_{N}^{+}$.

Finally, for a given state $y_{t}\in\mathcal{C}$, $t\in[0,r_{s}]$, and for the same integer $M$ previously adopted for the grid $\Omega_{M}^{-}$, we consider its discrete counterpart
\bd
\bs{y}_{t,M}=(y_{t}(\theta_{M,0}^{-}),\ldots,y_{t}(\theta_{M,M}^{-}))^{T}\in\mathcal{C}_{M}
\ed
with $\mathcal{C}_{M}:=\mathbb{C}^{M+1}$ the discrete counterpart of the state space $\mathcal{C}$. Correspondingly,
\bd
y_{t,M}(t):=\sum\limits_{j=0}^{M}\ell_{M,j}^{-}(t)y_{t}(\theta_{M,j}^{-}),\;t\in[-\tau,0],
\ed
is the relevant interpolating polynomial. It is not difficult to see that
\bd
y_{t,M}=\mathcal{L}_{M}^{-}y_{t}.
\ed
\br
Observe that $\mathcal{C}$ and $\mathcal{C}^{-}$ are the same space, as well as $\mathcal{C}_{M}$ and $\mathcal{C}_{M}^{-}$, all isomorphic to $\mathbb{C}^{M+1}$. However, we reserve to distinguish the notation for the relevant meaning and role. Conversely, $\mathcal{C}_{N}^{+}$ will be in general different from $\mathcal{C}_{M}^{-}$, since $N\neq M$ can be chosen, their role becoming clear after the analysis of convergence in Section \ref{s_convC}.
\er
\section{Discretization of the semigroup}\label{s_discrC}
We aim at finding a finite dimensional approximation $\bs{T}_{M,N}(r,s)$ of the evolution family $T(r,s)$ in (\ref{Tr}). We basically use {\it collocation}, together with {\it classic polynomial interpolation} as introduced in Section \ref{s_notC}. Briefly, in a discrete fashion we first transform the initial state into the solution on $[0,r_{s}]$ and, second, we restrict this latter when $r_{s}\geq\tau$, respectively prolong when $r_{s}<\tau$.

According to the notation set in Section \ref{s_notC} (but neglecting the reference to $r_{s}$ for simplicity), we first construct matrices $\bs{U}_{N,M}^{-}:\mathcal{C}_{M}^{-}\rightarrow\mathcal{C}_{N}^{+}$ and $\bs{U}_{N,N}^{+}:\mathcal{C}_{N}^{+}\rightarrow\mathcal{C}_{N}^{+}$ such that
\be\label{UNC}
\bs{U}_{N,N}^{+}\bs{p}_{N}^{+}=\bs{U}_{N,M}^{-}\bs{\varphi}_{M}^{-}
\ee
where $p_{N,M}\in\mathcal{C}^{\pm}$ is divided into
\be\label{pN+pN-}
\left\{
\setlength\arraycolsep{0.1em}\ba{ll}
\ds p_{N}^{+}(t):=p_{N,M}(t),&\;t\in[0,r_{s}]\\
\ds p_{M}^{-}(t):=p_{N,M}(t)=\varphi_{M}^{-}(t),&\;t\in[-\tau,0].
\ea
\right.
\ee
with $p_{N}^{+}$ determined by collocation of (\ref{IVPG_{s}yC}) on $\Omega_{N}^{+}$ with initial function $\varphi_{M}^{-}$:
\be\label{IVPG_{s}pNMC}
\left\{
\setlength\arraycolsep{0.1em}\ba{ll}
\ds (p_{N}^{+})'(\theta_{N,i}^{+})=(G_{s}p_{N,M})(\theta_{N,i}^{+}),&\;i=1,\ldots,N,\\
\ds p_{N}^{+}(0)=\varphi_{M}^{-}(0).
\ea
\right.
\ee
It is not difficult to check that
the above matrices have entries, respectively,
\bd
[\bs{U}_{N,N}^{+}]_{ij}:=
\left\{
\setlength\arraycolsep{2pt}\ba{ll}

\setlength\arraycolsep{2pt}\ba{l}
\ds1
\ea&
\left\{
\setlength\arraycolsep{2pt}\ba{l}
\ds i=0\\
\ds j=0
\ea
\right.\\

\setlength\arraycolsep{2pt}\ba{l}
\ds0
\ea&
\left\{
\setlength\arraycolsep{2pt}\ba{l}
\ds i=0\\
\ds j=1,\ldots,N
\ea
\right.\\

\left.\setlength\arraycolsep{2pt}\ba{l}
\ds(\ell_{N,j}^{+})'(\theta_{N,i}^{+})-a_{s}(\theta_{N,i}^{+})\delta_{ij}\\
\ds\quad-\int\limits_{-\theta_{N,i}^{+}}^{0}c_{s}(\theta_{N,i}^{+},\theta)\ell_{N,j}^{+}(\theta_{N,i}^{+}+\theta)d\theta
\ea\right\}&
\left\{
\setlength\arraycolsep{2pt}\ba{l}
\ds i=1,\ldots,N^{+}\\
\ds j=0,\ldots,N
\ea
\right.\\

\left.\setlength\arraycolsep{2pt}\ba{l}
\ds(\ell_{N,j}^{+})'(\theta_{N,i}^{+})-a_{s}(\theta_{N,i}^{+})\delta_{ij}\\
\ds\quad-b_{s}(\theta_{N,i}^{+})\ell_{N,j}^{+}(\theta_{N,i}^{+}-\tau)\\
\ds\quad-\int\limits_{-\tau}^{0}c_{s}(\theta_{N,i}^{+},\theta)\ell_{N,j}^{+}(\theta_{N,i}^{+}+\theta)d\theta
\ea\right\}&
\left\{
\setlength\arraycolsep{2pt}\ba{l}
\ds i=N^{+}+1,\ldots,N\\
\ds j=0,\ldots,N
\ea
\right.\\

\ea
\right.
\ed
and
\bd
[\bs{U}_{N,M}^{-}]_{ij}:=
\left\{
\setlength\arraycolsep{2pt}\ba{ll}

\setlength\arraycolsep{2pt}\ba{l}
\ds1
\ea&
\left\{
\setlength\arraycolsep{2pt}\ba{l}
\ds i=0\\
\ds j=0
\ea
\right.\\

\setlength\arraycolsep{2pt}\ba{l}
\ds0
\ea&
\left\{
\setlength\arraycolsep{2pt}\ba{l}
\ds i=0\\
\ds j=1,\ldots,M
\ea
\right.\\

\left.\setlength\arraycolsep{2pt}\ba{l}
\ds b_{s}(\theta_{N,i}^{+})\ell_{M,j}^{-}(\theta_{N,i}^{+}-\tau)\\
\ds\quad+\int\limits_{-\tau}^{-\theta_{N,i}^{+}}c_{s}(\theta_{N,i}^{+},\theta)\ell_{M,j}^{-}(\theta_{N,i}^{+}+\theta)d\theta
\ea\right\}&
\left\{
\setlength\arraycolsep{2pt}\ba{l}
\ds i=1,\ldots,N^{+}\\
\ds j=0,\ldots,M
\ea
\right.\\

\setlength\arraycolsep{2pt}\ba{l}
\ds0
\ea&
\left\{
\setlength\arraycolsep{2pt}\ba{l}
\ds i=N^{+}+1,\ldots,N\\
\ds j=0,\ldots,M,
\ea
\right.\\

\ea
\right.
\ed
where
\bd
N^{+}=N^{+}(r_{s},\tau):=\max\limits_{j=1,\ldots,N}\{\theta_{N,j}^{+}-\tau\leq0\}
\ed
and $\delta_{ij}$ is the Kronecker's delta.

Second, and independently of the model coefficients $a$, $b$ and $c$, we construct matrices $\bs{V}_{M,N}^{+}:\mathcal{C}_{N}^{+}\rightarrow\mathcal{C}_{M}$ and $\bs{V}_{M,M}^{-}:\mathcal{C}_{M}^{-}\rightarrow\mathcal{C}_{M}$ such that
\be\label{VNC}
\bs{y}_{r,M}=\bs{V}_{M,N}^{+}\bs{p}_{N}^{+}+\bs{V}_{M,M}^{-}\bs{\varphi}_{M}^{-}
\ee
by restriction of $p_{N,M}$ to $[r_{s}-\tau,r_{s}]$ when $r_{s}\geq\tau$, respectively prolongation by $\varphi_{M}$ when $r_{s}<\tau$. In particular, it is sufficient to define the above matrices with entries, respectively,
\bd
[\bs{V}_{M,N}^{+}]_{ij}:=
\left\{
\setlength\arraycolsep{2pt}\ba{ll}

\setlength\arraycolsep{2pt}\ba{l}
\ds\ell_{N,j}^{+}(r_{s}+\theta_{M,i}^{-})
\ea&
\left\{
\setlength\arraycolsep{2pt}\ba{l}
\ds i=0,\ldots,M^{-}\\
\ds j=0,\ldots,N
\ea
\right.\\

\setlength\arraycolsep{2pt}\ba{l}
\ds0
\ea&
\left\{
\setlength\arraycolsep{2pt}\ba{l}
\ds i=M^{-}+1,\ldots,M\\
\ds j=0,\ldots,N
\ea
\right.\\

\ea
\right.
\ed
and
\bd
[\bs{V}_{M,M}^{-}]_{ij}:=
\left\{
\setlength\arraycolsep{2pt}\ba{ll}

\setlength\arraycolsep{2pt}\ba{l}
\ds0
\ea&
\left\{
\setlength\arraycolsep{2pt}\ba{l}
\ds i=0,\ldots,M^{-}\\
\ds j=0,\ldots,M
\ea
\right.\\

\setlength\arraycolsep{2pt}\ba{l}
\ds\ell_{M,j}^{-}(r_{s}+\theta_{M,i}^{-})
\ea&
\left\{
\setlength\arraycolsep{2pt}\ba{l}
\ds i=M^{-}+1,\ldots,M\\
\ds j=0,\ldots,M,
\ea
\right.\\

\ea
\right.
\ed
where
\bd
M^{-}=M^{-}(r_{s},\tau):=\max\limits_{j=0,\ldots,M}\{r_{s}+\theta_{M,j}^{-}\geq0\},
\ed
with the convention that $\bs{V}_{M,N}^{+}$ is full and $\bs{V}_{M,M}^{-}$ is empty when $M^{-}=M$, i.e. for $r_{s}\geq\tau$.

Eventually, by setting $\bs{y}_{0,M}=\bs{\varphi}_{M}^{-}$, it follows from (\ref{UNC}) and (\ref{VNC}) that
\be\label{vecTNrC}
\bs{y}_{r_{s},M}=\bs{T}_{M,N}(r,s)\bs{y}_{0,M}
\ee
is the sought discrete approximation of (\ref{Tr}) with $\bs{T}_{M,N}(r,s):\mathcal{C}_{M}\rightarrow\mathcal{C}_{M}$ given by
\bd
\bs{T}_{M,N}(r,s)=\bs{V}_{M,N}^{+}(\bs{U}_{N,N}^{+})^{-1}\bs{U}_{N,M}^{-}+\bs{V}_{M,M}^{-}.
\ed
Standard approximation arguments ensure that $\bs{U}_{N,N}^{+}$ is invertible for sufficiently large $N$.

Aim of this research is to show how and under which conditions the (computable) spectrum of $\bs{T}_{M,N}(r,s)$ approximates that of $T(r,s)$.
\br
Let us observe that suitable quadrature rules have to be applied whenever the integrals in the above matrices $\bs{U}_{N,N}^{+}$ and $\bs{U}_{N,M}^{-}$ cannot be computed exactly. Of course, there is no choice but that of uniform distribution in order to exploit the nodal values already at disposal due to the varying integration windows $[\theta_{N,i}^{+}-\tau,\theta_{N,i}^{+}]$, $i=1,\ldots,N$. However, we feel like to advice the use of gaussian-type nodes (e.g. Gauss-Legendre or Chebyshev): at the price of extra technicalities in forming the matrix coefficients, this choice does not corrupt the spectral convergence that can be performed for the overall approximation as demonstrated in Section \ref{s_convC} (see \cite{bresisc05} for further details).
\er
\section{Convergence analysis}\label{s_convC}
The evolution family $T(r,s)$ in (\ref{Tr}) is an infinite dimensional map $T(r,s):\mathcal{C}\rightarrow\mathcal{C}$, contrary to its matrix discretization $\bs{T}_{M,N}(r,s):\mathcal{C}_{M}\rightarrow\mathcal{C}_{M}$ in (\ref{vecTNrC}). For comparison, it is therefore necessary to introduce an intermediate infinite dimensional, possibly finite rank, map $T_{N}(r,s):\mathcal{C}\rightarrow\mathcal{C}$. Set then
\be\label{TNr}
T_{N}(r,s)\varphi=(q_{N})_{r_{s}}
\ee
where $q_{N}\in\mathcal{C}^{\pm}$ is divided into
\be\label{qN+qN-}
\left\{
\setlength\arraycolsep{0.1em}\ba{ll}
\ds q_{N}^{+}(t):=q_{N}(t),&\;t\in[0,r_{s}]\\
\ds q_{N}^{-}(t):=q_{N}(t)=\varphi(t),&\;t\in[-\tau,0],
\ea
\right.
\ee
with $q_{N}^{+}$ determined by collocation of (\ref{IVPG_{s}yC}) on $\Omega_{N}^{+}$ with initial function $\varphi$:
\be\label{IVPG_{s}qNC}
\left\{
\setlength\arraycolsep{0.1em}\ba{ll}
\ds (q_{N}^{+})'(\theta_{N,i}^{+})=(G_{s}q_{N})(\theta_{N,i}^{+}),&\;i=1,\ldots,N,\\
\ds q_{N}^{+}(0)=\varphi(0).
\ea
\right.
\ee
Note that, in general, the collocation polynomials $q_{N}$ above and $p_{N,M}$ in (\ref{pN+pN-}) and (\ref{IVPG_{s}pNMC}) are different since relevant to different initial functions, $\varphi$ and $\varphi_{M}^{-}$, respectively. They coincide only when $\varphi\in\Pi_{M}^{-}$.

This Section is devoted first to provide in Section \ref{s_conv1} error bounds for the remainder $T(r,s)-T_{N}(r,s)$ in a suitable state space. Such errors will be measured in a {\it pointwise} sense in general (i.e. as applied to {\it a given} function in the chosen space), reserving to comment on the convergence in {\it norm} (i.e. as applied to {\it all} functions in the chosen space).

Second, in Section \ref{s_conv2} a link for the spectral elements of $\bs{T}_{M,N}(r,s)$ and $T_{N}(r,s)$ is studied, based on the relation
\be\label{TMNrTNr}
\bs{T}_{M,N}(r,s)=\mathcal{R}_{M}T_{N}(r,s)\mathcal{P}_{M},
\ee
not difficult to be verified.

Eventually, according to the theory developed in \cite{chat83}, it will be proved in Section \ref{s_conv3} that pointwise convergence of $T_{N}(r,s)$ to $T(r,s)$ in a suitable Banach space is a {\it mandatory} requirement for the approximation of the spectral elements of $T(r,s)$ by those of $T_{N}(r,s)$, a finite number of which eventually coincide with those of $\bs{T}_{M,N}(r,s)$ as it will be demonstrated by virtue of (\ref{TMNrTNr}). Instead, if the norm convergence is available (a much stringent requirement, rather difficult to happen), both the theories in \cite{chat83} and \cite{ggk90} can be applied. These comment hold for any numerical method that can be potentially used to form the matrix approximation $\bs{T}_{M,N}(r,s)$ (e.g. Runge-Kutta or Linear Multistep based schemes, see e.g. \cite{brean206,engsjna02}).

First to proceed, let us introduce some notation. For $\mathcal{C}$ and $Y$ normed linear spaces, let $\mathcal{B}(\mathcal{C},Y)$ be the set of linear and bounded operators from $\mathcal{C}$ to $Y$. Far any $A\in\mathcal{B}(\mathcal{C},Y)$ let
\bd
\|A\|_{\mathcal{C}\rightarrow Y}:=\sup\limits_{\psi\in\mathcal{C}}\frac{\|A\psi\|_{Y}}{\|\psi\|_{\mathcal{C}}}.
\ed
If $\mathcal{C}=Y$, then we will denote simply $\|A\|_{\mathcal{C}}$ for the operator norm of $A\in\mathcal{B}(\mathcal{C}):=\mathcal{B}(\mathcal{C},\mathcal{C})$.

For the convergence analysis, it will be often necessary to ask for more regularity than what so far demanded. To this aim, absolute continuity (see Definition \ref{d_AC} in Appendix \ref{s_app}) will be used and the space
\bd
\mathcal{C}_{A}:=W^{1,1}=\{\psi\in L\ :\ \psi'\in L\}
\ed
will be considered where
\bd
L:=L^{1}(-\tau,0;\mathbb{C}).
\ed
Spaces $\mathcal{C}_{A}^{-}$, $\mathcal{C}_{A}^{+}$ and $\mathcal{C}_{A}^{\pm}$ as well as $L^{-}$, $L^{+}$ and $L^{\pm}$ are similarly defined whether required. We will also resort to Lipschitz continuity, denoting the relevant spaces with ${\rm Lip}_{K}$ and ${\rm Lip}_{K}^{-}$, ${\rm Lip}_{K}^{+}$ and ${\rm Lip}_{K}^{\pm}$ whether required ($K$ denotes the Lipschitz constant).

For $\mathcal{C}$ we use $\|\cdot\|_{\mathcal{C}}$ as defined in (\ref{normC}), while for $\mathcal{C}_{A}$ we use
\be\label{normAC}
\|\psi\|_{\mathcal{C}_{A}}:=\|\psi\|_{L}+\|\psi'\|_{L}
\ee
with
\bd
\|\psi\|_{L}=\int\limits_{-\tau}^{0}|\psi(\theta)|d\theta.
\ed
Similar norms for $\mathcal{C}^{-}$, $\mathcal{C}^{+}$ and $\mathcal{C}^{\pm}$, $\mathcal{C}_{A}^{-}$, $\mathcal{C}_{A}^{+}$ and $\mathcal{C}_{A}^{\pm}$ and $L^{-}$, $L^{+}$ and $L^{\pm}$ are used whether required. With the above choices, all the spaces are of Banach type.

A series of technical and preparatory Lemmas are stated and proved in Appendix \ref{s_app}. Here we fully address only the main results. The constants $C$ appearing in the various statements (included those in Appendix \ref{s_app}), although different, always depend on $r$, $s$ and $\tau$ as well as on the coefficients $a$, $b$ and $c$.
\subsection{Convergence of $T_{N}(r,s)$ to $T(r,s)$}\label{s_conv1}
In this Section we study the convergence of the collocation under suitable hypotheses on the coefficients $a$, $b$ and $c$ and on the initial function.
\bas\label{cheb+}
Assume the nodes in $\Omega_{N}^{+}$ to be the zeros of the $N$-degree Chebyshev polynomial of the first kind in $(0,r_{s})$ \cite{davis75}, i.e.
\bd
\theta_{N,i}^{+}:=\frac{r_{s}}{2}\left(1-\cos{\left(\frac{(2i-1)\pi}{2N}\right)}\right),\;i=1,\ldots,N.
\ed
\eas


\bas\label{coef2}
Assume $a_{s}\in{\rm Lip}_{K_{a}}^{+}$, $b_{s}\in{\rm Lip}_{K_{b}}^{+}$ and $c_{s}(\cdot,\theta)\in{\rm Lip}_{K_{c}(\theta)}^{+}$ for all $\theta\in[-\tau,0]$ with $K_{c}\in L$. Assume, moreover, $\|c_{s}(\cdot,\theta)\|_{\mathcal{C}^{+}}\in L$ and denote
\bd
\|c_{s}\|_{\mathcal{C}^{+},L}:=\int\limits_{-\tau}^{0}\|c_{s}(\cdot,\theta)\|_{\mathcal{C}^{+}}d\theta.
\ed
\eas

\bth\label{t_xqN}
Let $G_{s}$ and $\mathcal{L}_{N}^{+}$ be defined by (\ref{G_{s}}) and (\ref{LN+}), respectively, and let $y\in\mathcal{C}^{\pm}$ be the solution of (\ref{IVPG_{s}yC}) with $\varphi\in\mathcal{C}$. Then, under Assumptions \ref{cheb+} and \ref{coef2} and for sufficiently large $N$, there exists a unique collocation solution $q_{N}\in\mathcal{C}^{\pm}$ given by (\ref{qN+qN-}) and (\ref{IVPG_{s}qNC}). Moreover,
\be\label{yqN}
\|y-q_{N}\|_{\mathcal{C}^{\pm}}\leq C\|\rho_{N}\|_{\mathcal{C^{+}}}
\ee
and
\be\label{yqN'}
\|y'-q_{N}'\|_{\mathcal{C}^{+}}\leq C'\|\rho_{N}\|_{\mathcal{C^{+}}}
\ee
hold where $C$ and $C'$ are constants independent of $N$ and $\varphi$ and
\be\label{rN}
\rho_{N}:=(I-\mathcal{L}_{N}^{+})y'\in\mathcal{C}^{+}.
\ee
\eth
\bpf
Given $\varphi\in\mathcal{C}$, let $y\in\mathcal{C}^{\pm}$ be the solution of (\ref{IVPG_{s}yC}) on $[-\tau,r_{s}]$. It is not difficult to see, by integration, that $y$ satisfies the functional equation in $\mathcal{C}^{\pm}$
\be\label{x}
y=u_{\varphi}+VG_{s}y
\ee
as soon as we consider $u\in\mathcal{C}^{\pm}$ as the function
\bd
u_{\varphi}(t):=\left\{
\setlength\arraycolsep{0.1em}\ba{ll}
\ds \varphi(0),&\;t\in[0,r_{s}]\\
\ds \varphi(t),&\;t\in[-\tau,0],
\ea
\right.
\ed
$V:\mathcal{C}^{+}\rightarrow\mathcal{C}^{\pm}$ as the integral operator
\be\label{V}
(Vy)(t):=\left\{
\setlength\arraycolsep{0.1em}\ba{ll}
\ds \int\limits_{0}^{t}y(\sigma)d\sigma,&\;t\in[0,r_{s}]\\
\ds 0,&\;t\in[-\tau,0]
\ea
\right.
\ee
and $G_{s}$ as given by (\ref{G_{s}}).

As for the collocation polynomial determined by (\ref{qN+qN-}) and (\ref{IVPG_{s}qNC}), observe that
\bd
\setlength\arraycolsep{0.1em}\ba{rcl}
q_{N}(t)&=&\ds q_{N}(0)+\int\limits_{0}^{t}(q_{N}^{+})'(\sigma)d\sigma\\
&=&\ds\varphi(0)+\int\limits_{0}^{t}\sum\limits_{j=1}^{N}m_{j}^{+}(\sigma)(q_{N}^{+})'(\theta_{N,j}^{+})d\sigma\\
&=&\ds\varphi(0)+\int\limits_{0}^{t}\sum\limits_{j=1}^{N}m_{j}^{+}(\sigma)\left[(G_{s}q_{N})(\theta_{N,j}^{+})\right]d\sigma,\;t\in[0,r_{s}],
\ea
\ed
holds where $m_{j}^{+}$'s are the Lagrange basis polynomials relevant to the nodes in $\Omega_{N}^{+}$ (different from the $\ell_{N,j}^{+}$'s). Then $q_{N}$ satisfies the functional equation in $\mathcal{C}^{\pm}$
\be\label{qN}
q_{N}=u_{\varphi}+V\mathcal{L}_{N}^{+}G_{s}q_{N}.
\ee

Consequently, the error function given by
\be\label{eN}
e_{N}:=y-q_{N}
\ee
(which is clearly zero in $[-\tau,0]$) satisfies the functional equation in $\mathcal{C}^{\pm}$
\be\label{eq1eN}
e_{N}=V\mathcal{L}_{N}^{+}G_{s}e_{N}+V\rho_{N},
\ee
as it can be seen by subtracting (\ref{qN}) from (\ref{x}) and by adding and subtracting $V\mathcal{L}_{N}^{+}G_{s}y$ in the result.

Now, it is not difficult to see that (\ref{eq1eN}) has a unique solution given by
\be\label{eq12}
e_{N}=V\overline{e}_{N}
\ee
where $\overline{e}_{N}\in\mathcal{C}^{+}$ is the unique solution of
the functional equation in $\mathcal{C}^{+}$
\be\label{eq2}
\overline{e}_{N}=\mathcal{L}_{N}^{+}G_{s}V\overline{e}_{N}+\rho_{N}
\ee
thanks to Lemma \ref{l_I-LG_{s}Vinv}. Moreover,
\bd
e_{N}'(t):=\left\{
\setlength\arraycolsep{0.1em}\ba{ll}
\ds \overline{e}_{N}(t),&\;t\in[0,r_{s}]\\
\ds 0,&\;t\in[-\tau,0].
\ea
\right.
\ed
The thesis is now straightforward.
\epf

Observe that Theorem \ref{t_xqN} shows that, under Assumptions \ref{cheb+} and \ref{coef2} and for sufficiently large $N$, the continuity of the initial function $\varphi$ is enough to guarantee the existence and uniqueness of the collocation solution as well as the error bound (\ref{yqN}). However, from (\ref{rN}) it is clear that continuity is not sufficient to ensure convergence, i.e.
\be\label{convxqN}
\lim\limits_{N\rightarrow\infty}\|y-q_{N}\|_{\mathcal{C}^{\pm}}=0,
\ee
for all $\varphi\in\mathcal{C}$. In fact, $G_{s}y$ involves also the initial function $\varphi$ due to the presence of the delay, hence it is only continuous and Faber's Theorem \cite{davis75,fab14} prevents convergence in all $\mathcal{C}^{\pm}$.

We now elaborate more, by resorting to the interpolation result in \cite{kry56}, which holds only under Assumption \ref{cheb+}.
\bpr\label{p_convrNAC}
Let $\rho_{N}$ be given by (\ref{rN}) under Assumption \ref{cheb+} and \ref{coef2}. If $\varphi\in\mathcal{C}_{A}$, then
\bd
\lim\limits_{N\rightarrow\infty}\|\rho_{N}\|_{\mathcal{C}^{+}}=0.
\ed
\epr
\bpf
If $\varphi\in\mathcal{C}_{A}$, then $G_{s}y\in\mathcal{C}_{A}^{\pm}$ follows easily from (\ref{G_{s}}). The thesis is then given by the result in \cite{kry56}, valid under Assumption \ref{cheb+}.
\epf

Proposition \ref{p_convrNAC} shows that absolute continuity is a minimal assumption for the initial function ensuring convergence as meant in (\ref{convxqN}).

Now we study the error for the approximated evolution family (\ref{TNr}) w.r.t. the exact one (\ref{Tr}).
\bco\label{c_TrTNr}
Let $T(r,s)$ be given by (\ref{Tr}). Then, under Assumptions \ref{cheb+} and \ref{coef2} and for sufficiently large $N$, $T_{N}(r,s)$ in (\ref{TNr}) is uniquely defined and, for any given $\varphi\in\mathcal{C}$,
\bd
\|(T(r,s)-T_{N}(r,s))\varphi\|_{\mathcal{C}}\leq C\|\rho_{N}\|_{\mathcal{C^{+}}}
\ed
holds where $\rho_{N}$ is given by (\ref{rN}) and $C$ is a constant independent of $N$ and $\varphi$. If, in addition, $\varphi\in\mathcal{C}_{A}$, then
\bd
\lim\limits_{N\rightarrow\infty}\|(T(r,s)-T_{N}(r,s))\varphi\|_{\mathcal{C}}=0.
\ed
\eco
\bpf
The thesis is straightforward by observing that
\bd
\|(T(r,s)-T_{N}(r,s))\varphi\|_{\mathcal{C}}\leq\|y-q_{N}\|_{\mathcal{C}^{\pm}}
\ed
holds for all $r\geq s$ and by applying Theorem \ref{t_xqN} and Proposition \ref{p_convrNAC}.
\epf

Corrollary \ref{c_TrTNr} ensures convergence for the approximated (still infinite dimensional) evolution family when applied to any absolutely continuous initial function in the topology induced by $\|\cdot\|_{\mathcal{C}}$. In order to apply the theory developed in \cite{chat83}, the basic requirement to be satisfied is that of {\it pointwise convergence in a Banach space}. Precisely (and according to Definition \ref{d_p} later on), given a Banach space $(B,\|\cdot\|_{B})$, a sequence of operators $\{A_{N}\}_{N=0}^{\infty}$ in $\mathcal{B}(B)$ is pointwise convergent to $A\in\mathcal{B}(B)$ if
\bd
\lim\limits_{N\rightarrow\infty}\|(A-A_{N})f\|_{B}=0
\ed
for all $f\in B$. The property is denoted by $A_{N}\xrightarrow{p}A$. Clearly, Theorem \ref{t_xqN} and the successive comment shows that $T_{N}(r,s)\xrightarrow{p}T(r,s)$ cannot hold in the natural state space $(\mathcal{C},\|\cdot\|_{\mathcal{C}})$. We now prove that it holds instead if we choose as the state space $(\mathcal{C}_{A},\|\cdot\|_{\mathcal{C}_{A}})$ with $\|\cdot\|_{\mathcal{C}_{A}}$ given by (\ref{normAC}).
\bth\label{t_TrTNrAC}
Let $T(r,s)$ and $T_{N}(r,s)$ be given by (\ref{Tr}) and (\ref{TNr}), respectively, under Assumptions \ref{cheb+} and \ref{coef2} and for sufficiently large $N$. Then $T_{N}(r,s)\xrightarrow{p}T(r,s)$ in $(\mathcal{C}_{A},\|\cdot\|_{\mathcal{C}_{A}})$. Moreover,
\bd
\sup\limits_{N\in\mathbb{N}}\|T_{N}(r,s)\|_{\mathcal{C}_{A}}<\infty.
\ed
\eth
\bpf
We use the same notation as in the proof of Theorem \ref{t_xqN}. Since
\bd
\|e_{N}\|_{\mathcal{C}_{A}}\leq r_{s}(\|e_{N}\|_{\mathcal{C}^{\pm}}+\|\overline{e}_{N}\|_{\mathcal{C}^{+}}),
\ed
pointwise convergence follows by virtue of Proposition \ref{p_convrNAC}. The last assertion follows from the Banach-Steinhaus Theorem \cite[Theorem 3.1]{chat83}.
\epf

Let us underline again that the sequence of operators $\{T_{N}(r,s)\}_{N=0}^{\infty}$ is made of infinite dimensional maps which have finite rank only when $r_{s}\geq\tau$. In fact, when $r_{s}<\tau$, the image $T_{N}(r,s)\varphi$ contains a piece of the initial function $\varphi$, precisely
\bd
T_{N}(r,s)\varphi(t)=\left\{
\setlength\arraycolsep{0.1em}\ba{ll}
\ds q_{N}(t),&\;t\in[0,r_{s}]\\
\ds \varphi(t),&\;t\in[r_{s}-\tau,0],
\ea
\right.
\ed
which, in general, is not a polynomial. However, for $r_{s}<\tau$, $T(r,s)$ is neither compact \cite{hvl93}. Therefore, for $r_{s}\geq\tau$ we have constructed a sequence of finite rank approximations to the exact evolution family, which can be proved to remain compact in $\mathcal{B}(\mathcal{C}_{A})$ as well by standard arguments.

We conclude the Section with an important remark on the convergence in norm.
\br\label{r_norm}
Absolute continuity is sufficient to provide pointwise convergence, but for norm convergence more regularity is needed. In fact, Proposition \ref{p_convrNAC} states that the interpolation remainder $\|\rho_{N}\|_{\mathcal{C}^{+}}$ vanishes, while for norm convergence
\bd
\|\rho_{N}\|_{\mathcal{C}^{+}}\leq g(N)\|y\|_{\mathcal{C}^{\pm}}
\ed
with $g(N)$ vanishing independently of $y$ would be necessary.

It is not difficult to see from the proof of Theorem \ref{t_TrTNrAC} and from standard interpolation results as applied to $\rho_{N}$ that, if we further restrict to the state space $\mathcal{C}^{1}:=C^{1}(-\tau,0;\mathbb{C})$ with its natural norm
\bd
\|\psi\|_{\mathcal{C}^{1}}:=\|\psi\|_{\mathcal{C}}+\|\psi'\|_{\mathcal{C}},
\ed
then we obtain convergence in norm, i.e.
\bd
\lim\limits_{N\rightarrow\infty}\|T(r,s)-T_{N}(r,s)\|_{\mathcal{C}^{1}}=0.
\ed
Further comments on this result will be made at the end of Section \ref{s_conv3}.
\er
\subsection{Relation between $T_{N}(r,s)$ and $\bs{T}_{M,N}(r,s)$}\label{s_conv2}
With an eye kept on the spectral elements, we analyze now the relation between the infinite dimensional approximation $T_{N}(r,s)$ defined in (\ref{TNr}) and the finite dimensional one $\bs{T}_{M,N}(r,s)$ introduced in Section \ref{s_discrC} through (\ref{vecTNrC}). To this aim, $N$ is considered fixed throughout the whole Section. The following result is fundamental to the scope.

\bpr\label{p_TNrTMNr}
Let the matrix $\bs{A}_{M}:\mathcal{C}_{M}\rightarrow\mathcal{C}_{M}$ be given and define the operator $A_{M}:=\mathcal{P}_{M}^{-}\bs{A}_{M}\mathcal{R}_{M}^{-}\in\mathcal{B}(\mathcal{C})$ through the restriction and prolongation maps $\mathcal{R}_{M}^{-}$ and $\mathcal{P}_{M}^{-}$, respectively, introduced in Section \ref{s_notC}. Then $\bs{A}_{M}$ and $A_{M}$ have the same nonzero eigenvalues with the same geometric and partial multiplicities.
\epr
\bpf
Let $\mu\in\mathbb{C}\setminus\left\{0\right\}$. We prove that
\bd
\mathcal{P}_{M}^{-}\ker(\mu\bs{I}_{M}-\bs{A}_{M})=\ker(\mu I-A_{M}).
\ed
Let $\bs{v}\in\mathcal{C}_{M}$, $\bs{v}\neq0$, such that 
\bd
\bs{A}_{M}\bs{v}=\mu\bs{v}. 
\ed
Let $\psi=\mathcal{P}_{M}^{-}\bs{v}$. Then
\bd
A_{M}\psi=\mathcal{P}_{M}^{-}\bs{A}_{M}\bs{v}=\mathcal{P}_{M}^{-}\mu\bs{v}=\mu\psi. 
\ed
Vice versa, let $\psi\in\mathcal{C}$, $\psi\neq0$, such that 
\bd
A_{M}\psi=\mu\psi.
\ed
Since $\mu\neq0$, we have
\bd
\psi=\mathcal{P}_{M}^{-}\bs{v} 
\ed
where 
\bd
\bs{v}=\frac{1}{\mu}\bs{A}_{M}\mathcal{R}_{M}^{-}\psi\in\mathcal{C}_{M}. 
\ed
Hence
\bd
\mathcal{P}_{M}^{-}\bs{A}_{M}\bs{v}=\mu \mathcal{P}_{M}^{-}\bs{v} 
\ed
and then
\bd
\bs{A}_{M}\bs{v}=\mu\bs{v}. 
\ed

Now, we prove that 
\bd
\dim\ker(\mu\bs{I}_{M}-\bs{A}_{M})=\dim\ker(\mu I-A_{M})
\ed
Let $\bs{v}_{1},\ldots,\bs{v}_{g}$ be linearly independent elements of $\ker(\mu\bs{I}_{M}-\bs{A}_{M})$. We have shown above that
\bd
\psi_{i}=\mathcal{P}_{M}^{-}\bs{v}_{i}\in\ker(\mu I-A_{M}),\ i=1,\ldots,g. 
\ed
Since
\bd
\sum\limits_{i=1}^{g}\alpha_{i}\psi_{i}=\mathcal{P}_{M}^{-}
\sum\limits_{i=1}^{g}\alpha_{i}\bs{v}_{i}=0\Longrightarrow
\sum\limits_{i=1}^{g}\alpha_{i}\bs{v}_{i}=0\Longrightarrow\alpha_{i}=0,\ i=1,\ldots,g,
\ed
the elements $\psi_{1},\ldots,\psi_{g}$ are linearly independent. Viceversa, let $\psi_{1},\ldots,\psi_{g}$ be linearly independent elements of $\ker(\mu I-A_{M})$. We have shown above that
\bd
\psi_{i}=\mathcal{P}_{M}^{-}\bs{v}_{i},\ i=1,\ldots,g, 
\ed
where $\bs{v}_{i}\in\ker(\mu\bs{I}_{M}-\bs{A}_{M})$. Since
\bd
\sum\limits_{i=1}^{g}\alpha_{i}\bs{v}_{i}=0\Longrightarrow \mathcal{P}_{M}^{-}\sum\limits_{i=1}^{g}\alpha_{i}\bs{v}_{i}=\sum\limits_{i=1}^{g}\alpha_{i}\psi_{i}\Longrightarrow\alpha_{i}=0,\ i=1,\ldots,g,
\ed
the elements $\bs{v}_{1},\ldots,\bs{v}_{g}$ are linearly independent.

Finally, let us prove that there is a one-to-one correspondence between Jordan chains of $\bs{A}_{M}$ and $A_{M}$. Let $\bs{v}_{1},\ldots,\bs{v}_{g}$ be a Jordan chain of $\bs{A}_{M}$. Then $\psi_{1}=\mathcal{P}_{M}^{-}\bs{v}_{1},\ldots,\psi_{g}=\mathcal{P}_{M}^{-}\bs{v}_{g}$ is a Jordan chain for $A_{M}$. In fact, $\psi_{1}$ is an eigenvector of $A_{M}$ and
\bd
{\setlength\arraycolsep{2pt}\ba{rcl}
(\mu I-A_{M})\psi_{i+1}&=&(\mu I-A_{M})\mathcal{P}_{M}^{-}\bs{v}_{i+1}\\
&=&\mathcal{P}_{M}^{-}(\mu\bs{I}_{M}-\bs{A}_{M})\bs{v}_{i+1}\\
&=&\bs{v}_{i},\ i=0,\ldots,g-1.
\ea}
\ed
Vice versa, let $\psi_{1},\ldots,\psi_{g}$ be a Jordan chain of $
A_{M}$. We have seen that 
\bd
\psi_{1}=\mathcal{P}_{M}^{-}\bs{v}_{1} 
\ed
for some $\bs{v}_{1}\in\mathcal{C}_{N}$ eigenvector of $\bs{A}_{M}$ (since $\psi_{1} $ is an eigenvector of $A_{M}$). Note that, for $i=0,\ldots,g-1$, if $\psi_{i}=\mathcal{P}_{M}^{-}\bs{v}_{i}$ for some $\bs{v}_{i}\in\mathcal{C}_{M}$, then $\psi_{i+1}=\mathcal{P}_{M}^{-}\bs{v}_{i+1}$ for some $\bs{v}_{i+1}\in\mathcal{C}_{M}$. In fact, 
\bd
(\mu I-A_{M})\psi_{i+1}=\psi_{i}=\mathcal{P}_{M}^{-}\bs{v}_{i} 
\ed
and so
\bd
\psi_{i+1}=\mathcal{P}_{M}^{-}\left(\frac{1}{\mu}(\bs{A}_{M}\mathcal{R}_{M}^{-}\psi_{i+1}-\bs{v}_{i})\right). 
\ed
We conclude that 
\bd
\psi_{i}=\mathcal{P}_{M}^{-}\bs{v}_{i},\ i=1,\ldots,g,
\ed
for some $\bs{v}_{i}\in\mathcal{C}_{M}$, i.e. $\bs{v}_{1},\ldots,\bs{v}_{g}$ is a Jordan chain for $\bs{A}_{M}$. In fact,
\bd
{\setlength\arraycolsep{2pt}\ba{rcl}
\mathcal{P}_{M}^{-}\bs{v}_{i} =\psi_{i}&=&(\mu I-A_{M})\psi_{i+1}\\
&=&(\mu I-A_{M})\mathcal{P}_{M}^{-}\bs{v}_{i+1}\\
&=&\mathcal{P}_{M}^{-}(\mu\bs{I}_{M}-\bs{A}_{M})\bs{v}_{i+1},\ i=0,\ldots,g-1,
\ea}
\ed
and then
\bd
\bs{v}_{i}=(\mu\bs{I}_{M}-\bs{A}_{M})\bs{v}_{i+1},\ i=0,\ldots,g-1. 
\ed
\epf

Proposition \ref{p_TNrTMNr} as applied to $\bs{T}_{M,N}(r,s)$ shows that this latter and the operator $\mathcal{P}_{M}^{-}\bs{T}_{M,N}(r,s)\mathcal{R}_{M}^{-}$ have the same nonzero eigenvalues and multiplicities. Moreover, in the proof it can be seen how the eigenvectors of $\bs{T}_{M,N}(r,s)$ and eigenfunctions of $\mathcal{P}_{M}^{-}\bs{T}_{M,N}(r,s)\mathcal{R}_{M}^{-}$ are related to each other via the restriction and prolongation operators.

Now, observe that (\ref{TMNrTNr}) and (\ref{LM-}) give
\bd
\mathcal{P}_{M}^{-}\bs{T}_{M,N}(r,s)\mathcal{R}_{M}^{-}=\mathcal{L}_{M}^{-}T_{N}(r,s)\mathcal{L}_{M}^{-}.
\ed
In view of this and of the results of Section \ref{s_conv1}, we need only to investigate the relation between the spectrum of $\mathcal{L}_{M}^{-}T_{N}(r,s)\mathcal{L}_{M}^{-}$ and that of $T_{N}(r,s)$.

\bth\label{t_LTNrLTNrspec}
Let $T_{N}(r,s)$ be defined by (\ref{TNr}) under Assumptions \ref{cheb+} and \ref{coef2} and with $r_{s}\geq\tau$ and let $\mathcal{L}_{M}^{-}$ be defined by (\ref{LM-}) with $M\geq N$. Then the spectral elements of $T_{N}(r,s)$ coincide with those of $\mathcal{L}_{M}^{-}T_{N}(r,s)\mathcal{L}_{M}^{-}$.
\eth
\bpf
On the one hand, we soon observe that for $M\geq N$ and $r_{s}\geq\tau$ we have
\bd
\mathcal{L}_{M}^{-}T_{N}(r,s)=T_{N}(r,s).
\ed
In fact, for $r_{s}\geq\tau$ the range of $T_{N}(r,s)$ is $\Pi_{N}^{-}$ and $\mathcal{L}_{M}^{-}\Pi_{N}^{-}=\Pi_{N}^{-}$ for $M\geq N$.

On the other hand, let $\psi$ be an eigenfunction of $T_{N}(r,s)$ associated to the eigenvalue $\mu$, i.e.
\bd
T_{N}(r,s)\psi=\mu\psi.
\ed
Then, for $r_{s}\geq\tau$, it must be $\psi\in\Pi_{N}^{-}$. Consequently,
\bd
T_{N}(r,s)\mathcal{L}_{M}^{-}\psi=T_{N}(r,s)\psi
\ed
whenever $M\geq N$. Arguments similar to those used to prove Proposition \ref{p_TNrTMNr} complete the proof.
\epf

As an immediate consequence of Theorem \ref{t_LTNrLTNrspec}, it is sufficient to take $M=N$ in the construction of the matrix $\bs{T}_{M,N}(r,s)$ in Section \ref{s_discrC}, so to keep the computational effort as low as possible. Moreover, it is interesting to observe that no special assumption on the distribution of the nodes $\Omega_{M}^{-}$ has to be made, contrary to Assumption \ref{cheb+} for $\Omega_{N}^{+}$. In fact, these nodes serve merely for representing $T_{N}(r,s)$ in finite dimension, and not for approximation reasons as those in $\Omega_{N}^{+}$. However, Theorem \ref{t_LTNrLTNrspec} holds only for $r_{s}\geq\tau$, which is the case of interest since, as already remarked, this condition guarantees the compactness of $T(r,s)$, and hence the nice and known properties of its spectrum. Nevertheless, it is not difficult to see that pointwise convergence $\mathcal{L}_{M}^{-}T_{N}(r,s)\mathcal{L}_{M}^{-}\xrightarrow{p}T_{N}(r,s)$ holds in all $(\mathcal{C}_{A},\|\cdot\|_{\mathcal{C}_{A}})$ as $M\rightarrow\infty$ and $N$ fixed as long as $r_{s}\geq\tau$.
\subsection{Computation and convergence for the spectrum of $T(r,s)$}\label{s_conv3}
From now on we assume $r_{s}\geq\tau$ and $M\geq N$.

In Section \ref{s_conv2} it is shown that, for a fixed $N$, the spectral elements of $T_{N}(r,s)$ coincide with those of the matrix $\bs{T}_{M,N}(r,s)$, through which they can be effectively and efficiently computed as the result of a standard algebraic eigenvalue problem in finite dimension. The last step consists in proving that these elements converge to those of $T(r,s)$ in the limit as $N\rightarrow\infty$. To this aim we apply the theory developed in \cite{chat83}, that we briefly recall for what concerns its basic facts and as adapted to the notation of the present manuscript.

Following \cite{chat83}, let $\mu$ be a nonzero and isolated eigenvalue of $T(r,s)$ with finite algebraic multiplicity $m$, geometric multiplicity $g$ and ascent $\ell$. Let $\Delta$ be a neighborhood of $\mu$ such that $\mu$ is the only eigenvalue of $T(r,s)$ in it. Let $\Gamma$ be a closed Jordan curve isolating $\mu$ and drawn in $\Delta$. Observe that all this makes sense since, for $r_{s}\geq\tau$, $T(r,s)$ is compact and hence it has only point spectrum with nontrivial isolated eigenvalues of finite algebraic multiplicity plus possibly $0$ as accumulation point \cite[Theorem 2.34]{chat83}, \cite[Chapter 7]{hvl93}. To such a $\mu$ it is associated the spectral projection
\bd
P:=\frac{1}{2\pi{\rm i}}\oint_{\Gamma}(T(r,s)-zI)^{-1}dz
\ed
and the relevant generalized eigenspace $\mathcal{M}:=P\mathcal{C}_{A}$. Recall that $m=\dim{M}$, $g=\dim{{\rm ker}(T(r,s)-\mu I)}$, while $(T(r,s)-\mu I)P$ is nilpotent and $\ell$ is the maximum integer such that $((T(r,s)-\mu I)P)^{\ell}=0$. Moreover, $1\leq g,\ell\leq m<\infty$. If $Q$ is another porjection, then for $\mathcal{N}:=Q\mathcal{C}_{A}$ it is defined the gap
\bd
\Theta(\mathcal{M},\mathcal{N}):=\max\{{\rm dist}(\mathcal{M},\mathcal{N}),{\rm dist}(\mathcal{N},\mathcal{M})\}
\ed
with
\bd
{\rm dist}(\mathcal{M},\mathcal{N}):=\sup\limits_{\substack{\psi\in\mathcal{M}\\ \|\psi\|_{\mathcal{C}_{A}}=1}}{\rm dist}(\psi,\mathcal{N})\neq\sup\limits_{\substack{\xi\in\mathcal{N}\\ \|\xi\|_{\mathcal{C}_{A}}=1}}{\rm dist}(\xi,\mathcal{M})=:{\rm dist}(\mathcal{N},\mathcal{M}).
\ed

The following definitions are necessary, also to understand the role of the analysis so far conducted.

\bdf[pointwise convergence \cite{chat83}]\label{d_p}
$T_{N}(r,s)\xrightarrow{p}T(r,s)$ iff, for all $\varphi\in\mathcal{C}_{A}$, $T_{N}(r,s)\varphi\rightarrow T(r,s)\varphi$ as $N\rightarrow\infty$.
\edf

\bdf[stable convergence \cite{chat83}]\label{d_s}
$T_{N}(r,s)-zI\xrightarrow{s}T(r,s)-zI$ for $z\in\Delta\setminus\{\mu\}$ iff
\bi
\item[(i)]$T_{N}(r,s)\xrightarrow{p}T(r,s)$;
\item[(ii)]$\exists K>0$, $\exists\overline{N}$ s.t. for $N>\overline{N}$, $(T_{N}(r,s)-zI)^{-1}\in\mathcal{B}(\mathcal{C}_{A})$ and $\|(T_{N}(r,s)-zI)^{-1}\|_{\mathcal{C}_{A}}\leq K$.
\ei
\edf

\bdf[strong stable convergence \cite{chat83}]\label{d_ss}
$T_{N}(r,s)-zI\xrightarrow{ss}T(r,s)-zI$ for $z\in\Delta\setminus\{\mu\}$ iff
\bi
\item[(i)]$T_{N}(r,s)-zI\xrightarrow{s}T(r,s)-zI$ for $z\in\Delta\setminus\{\mu\}$;
\item[(ii)]$\dim{P_{N}\mathcal{C}_{A}}=m$ for $N$ large enough.
\ei
If $T_{N}(r,s)-zI\xrightarrow{ss}T(r,s)-zI$ in $\Delta$, then $T_{N}(r,s)$ is said a strongly stable approximation of $T(r,s)$ in $\Delta$.
\edf

Definitions \ref{d_s} and \ref{d_ss} above hold similarly for all $z\in\Gamma$ if both conditions $(i)$ hold for $z\in\Gamma$.

A first important result follows.

\bpr[Proposition 5.6 in \cite{chat83}]\label{p_chat5.6}
If $T_{N}(r,s)$ is a strongly stable approximation of $T(r,s)$ in $\Delta$, then, $T_{N}(r,s)$ has in $\Delta$, for $N$ large enough, exactly $m$ eigenvalues, counting their multiplicities.
\epr

Let us call $\mu_{jN}$, $j=1,\ldots,m$, such eigenvalues, let $\mu_{N}$ represent any among the distinct ones and let $\psi_{N}$, $\|\psi_{N}\|_{\mathcal{C}_{A}}=1$, be an associated eigenfunction. For this $\mu_{N}$, let $P_{N}$ and $\mathcal{M}_{N}$ be, respectively, the relevant spectral projection and generalized eigenspace as previously introduced. Since $T(r,s)$ is not self-adjoint in general, $\mu$ is best approximated by the arithmetic mean
\bd
\hat{\mu}_{N}:=\sum\limits_{j=1}^{m}\mu_{jN}.
\ed

We recall now the two fundamental results, which involve the quantity
\be\label{epsN}
\varepsilon_{N}:=\|(T(r,s)-T_{N}(r,s))P\|_{\mathcal{C}_{A}}.
\ee
Observe that $\varepsilon_{N}$ is the remainder of the approximation as restricted to the generalized eigenspace, rather than on the whole space.

\bth[Theorem 6.6 in \cite{chat83}]\label{t_chat6.6}
If $T_{N}(r,s)$ is a strongly stable approximation of $T(r,s)$ in $\Delta$, then, for $N$ large enough, the quantities
\bd
\setlength\arraycolsep{0.1em}\ba{ll}
\ds \|(I-P)\psi_{N}\|_{\mathcal{C}_{A}},&\;\textrm{for }\psi_{N}\in\mathcal{M}_{N},\;\|\psi_{N}\|_{\mathcal{C}_{A}}=1,\\
\ds \|(I-P_{N})\psi\|_{\mathcal{C}_{A}},&\;\textrm{for }\psi\in\mathcal{M},\\
\ds\Theta(\mathcal{M},\mathcal{M}_{N}),&\\
\ds\mu-\hat{\mu}_{N},&\\
\ds\frac{1}{\mu}-\frac{1}{m}\left(\sum\limits_{j=1}^{m}\frac{1}{\mu_{jN}}\right)&
\ea
\ed
are at least of order $\varepsilon_{N}$.
\eth

\bth[Theorem 6.7 in \cite{chat83}]\label{t_chat6.7}
If $T_{N}(r,s)$ is a strongly stable approximation of $T(r,s)$ in $\Gamma$, then, for $N$ large enough,
\bd
\setlength\arraycolsep{0.1em}\ba{ll}
\ds\max\limits_{j=1,\ldots,m}|\mu-\mu_{jN}|=O(\varepsilon_{N}^{1/\ell}),\\
\ds\min\limits_{j=1,\ldots,m}|\mu-\mu_{jN}|=O(\varepsilon_{N}^{g/m}),\\
\ds{\rm dist}(\psi_{N},{\rm ker}(T(r,s)-\mu I))=O(\varepsilon_{N}^{1/\ell}).
\ea
\ed
\eth

In order for the previous results to hold true, we need now to verify the (only) hypothesis of strongly stable convergence according to Definition \ref{d_ss}. This, in turn, requires $(ii)$ in Definition \ref{d_ss}, $(ii)$ in Definition \ref{d_s} and pointwise convergence.

The latter is provided by Theorem \ref{t_TrTNrAC}, justifying all the analysis performed in Section \ref{s_conv1}.

As for $(ii)$ in Definition \ref{d_s}, it is enough to observe that $z\in\Delta\setminus\{\mu\}$ implies $z\in\rho(T_{N}(r,s))$, the resolvent set of $T_{N}(r,s)$, hence $(ii)$ holds by the definition of this latter.

Eventually, as for $(ii)$ in Definition \ref{d_ss}, we first recall the following definitions and results about the spectral projections $P_{N}$ and $P$. Below we use $B:=\{\psi\in\mathcal{C}_{A}\ :\ \|\psi\|_{\mathcal{C}_{A}}\leq1\}$.
 
\bdf[collectively compact convergence \cite{chat83}]\label{d_cc}
$P_{N}\xrightarrow{cc}P$ iff
\bi
\item[(i)]$P_{N}\xrightarrow{p}P$;
\item[(ii)]the set $\bigcup\limits_{N=1}^{\infty}(P-P_{N})B$ is relatively compact in $\mathcal{C}_{A}$.
\ei
\edf

\bdf[compact convergence \cite{chat83}]\label{d_c}
$P_{N}\xrightarrow{c}P$ iff
\bi
\item[(i)]$P_{N}\xrightarrow{p}P$;
\item[(ii)]for any sequence $\{\xi_{N}\}_{N=1}^{\infty}$ in $B$, the sequence $\{(P-P_{N})\xi_{N}\}_{N=1}^{\infty}$ is relatively compact in $\mathcal{C}_{A}$.
\ei
\edf

\bth[Theorem 3.9 in \cite{kress89}]\label{t_kress3.9}
The projections $P$ and $P_{N}$ are compact.
\eth

\bpr[\cite{ans71}, Proposition 3.13 in \cite{chat83}]\label{p_chat3.13}
For projections $P$ and $P_{N}$ such that $P$ is compact, the following are equivalent:
\bi
\item[(i)]$P_{N}\xrightarrow{p}P$ and $\dim{P_{N}\mathcal{C}_{A}}=\dim{P\mathcal{C}_{A}}<\infty$ for $N$ large enough;
\item[(ii)]$P_{N}\xrightarrow{cc}P$.
\ei
\epr

Theorem \ref{t_kress3.9} ensures that $P_{N}-P$ is compact, hence Definitions \ref{d_cc} and \ref{d_c} are equivalent, \cite[p.125]{chat83}. Moreover, $(ii)$ in Definition \ref{d_c} is trivially satisfied. Then, by virtue of Proposition \ref{p_chat3.13}, $(ii)$ in Definition \ref{d_ss} is guaranteed as soon as $P_{N}\xrightarrow{p}P$. As for this latter, for any given $\varphi\in\mathcal{C}_{A}$,
\bd
\setlength\arraycolsep{0.1em}\ba{rcl}
\ds (P-P_{N})\varphi&=&\ds\frac{1}{2\pi{\rm i}}\oint_{\Gamma}[(T(r,s)-zI)^{-1}-(T_{N}(r,s)-zI)^{-1}]\varphi dz\\
&=&\ds\frac{1}{2\pi{\rm i}}\oint_{\Gamma}(T_{N}(r,s)-zI)^{-1}[T_{N}(r,s)-T(r,s)](T(r,s)-zI)^{-1}\varphi dz
\ea
\ed
and hence $P_{N}\xrightarrow{p}P$ follows from $T_{N}(r,s)\xrightarrow{p}T(r,s)$, again guaranteed by Theorem \ref{t_TrTNrAC}.

Now we are able to comment on the convergence rate. Thanks to Theorems \ref{t_chat6.6} and \ref{t_chat6.7}, the error between the spectral elements of $T_{N}(r,s)$ and $T(r,s)$ decreases as fast as (\ref{epsN}). Thanks to Theorem \ref{t_xqN}, this latter is governed by $\|\rho_{N}\|_{\mathcal{C}^{+}}$ with $\rho_{N}$ given by (\ref{rN}) for $y$ solution of (\ref{IVPG_{s}yC}) with a generalized eigenfunction associated to $\mu$ as initial function. As it is well-known, the eigenfunctions of $T(r,s)$ are analytic \cite{hvl93}. This, together with Jackson's type Theorems, ensures {\it spectral accuracy}, i.e. a convergence of infinite order:
\bd
\varepsilon_{N}\leq \left(\frac{K}{N}\right)^{N}
\ed
where $K=K(|\mu|)$ is a constant independent of $N$.

We resume all the spectral convergence analysis in the following, where $\Sigma_{0}(\cdot)$ denotes formally the spectral elements of an operator (i.e. any of the quantities mentioned in Theorems \ref{t_chat6.6} or \ref{t_chat6.7}).
\bth\label{t_final}
Let $T(r,s)$ be given by (\ref{Tr}) with $r_{s}\geq\tau$. Under Assumption \ref{cheb+} and for sufficiently large $N$, $T_{N}(r,s)$ in (\ref{TNr}) is uniquely defined. For $M\geq N$, set $\bs{T}_{M,N}(r,s)=\mathcal{R}_{M}T_{N}(r,s)\mathcal{P}_{M}$. Then
\bi
\item[(i)]for a fixed $N$, $\Sigma_{0}(\bs{T}_{M,N}(r,s))=\Sigma_{0}(T_{N}(r,s))$;
\item[(ii)]as $N\rightarrow\infty$, $\Sigma_{0}(T_{N}(r,s))\rightarrow\Sigma_{0}(T(r,s))$ with spectral accuracy.
\ei
\eth
Let us underline that the equality stated in $(i)$ above is true for the eigenvalues and relevant multiplicities, while it hides the restriction and prolongation operators as far as the eigenvectors and eigenfunctions are concerned.

\br
By recalling Remark \ref{r_norm}, the result about convergence in norm there stated may be used to apply the theory developed in \cite{ggk90} in order to obtain the similar results as given in this Section for the approximation of the spectral elements. Such a convergence analysis is similar to that developed in \cite{brenm09} for partial retarded functional differential equations, and it is not based on the theory developed in \cite{chat83}. However, it is worthy to underline that we chose to follow \cite{chat83} since it requires the least possible restriction of the natural state space, i.e. $\mathcal{C}_{A}$ rather than $\mathcal{C}^{1}$.
\er

\section{Applications}\label{s_appl}
It is well-known (e.g. \cite{hvl93,diekmann95}) that the asymptotic stability of stationary solutions can be characterized through the knowledge of the spectrum of $T(r,s)$. In particular, the following two situations are certainly worthy to be mentioned in view of their importance in a wide class of applications.
\bi
\item When (\ref{dde}) is the result of the linearization of a nonlinear system of DDEs around an equilibrium solution, then the coefficients are autonomous, i.e. for all $t\in[0,r]$, $a(t)=a$, $b(t)=b$ and $c(t,\theta)=c(\theta)$ for all $\theta\in[-\tau,0]$ and, consequently, the evolution family $\{T(r,s)\}_{r\geq s}$ reduces to the standard $C_{0}$-semigroup of solution operators $\{T(r)\}_{r\geq0}$. Moreover, for $r\geq\tau$, such an equilibria is asymptotically stable iff $|\mu|<1$ for all the eigenvalues $\mu$ of $T(r)$, the so-called {\it multipliers}.
\item When (\ref{dde}) is the result of the linearization of a nonlinear system of DDEs around a limit cycle with period $\omega$, then the coefficients are periodic with the same period, i.e. for all $t\in[0,r]$, $a(t+\omega)=a(t)$, $b(t+\omega)=b(t)$ and $c(t+\omega,\theta)=c(t,\theta)$ for all $\theta\in[-\tau,0]$ and, consequently, the spectral properties of the evolution family $\{T(r,s)\}_{r\geq s}$ can be studied through the monodromy operator $U(\omega):=T(\omega,0)$. Moreover, for any $\omega\geq\tau$, such a periodic orbit is asymptotically stable iff $|\mu|<1$ for all the eigenvalues $\mu$ of $U(\omega)$, the so-called {\it Floquet multipliers}.
\ei
It is then clear that the results of the present work can be straightforwardly used in both cases without any additional effort. Let us observe, moreover, that for periodic problems where the period of the coefficients $\omega$ is less then the maximum delay $\tau$ the monodromy operator is not compact, even though a sufficiently large power $k$ of it is so, precisely $U^{k}(\omega)=U(k\omega)$ with some integer $k$ such that $k\omega\geq\tau$. Then, for stability purposes, spectral approximations of the latter are sufficient.
\appendix
\section{Appendix}\label{s_app}
\bdf\label{d_AC}
Let $(\mathcal{C},d)$ be a metric space. A function $\mathbb{R}\supset I\ni t\mapsto f(t)\in\mathcal{C}$ is {\rm absolutely continuous} if for any $\varepsilon>0$ there exists a $\delta=\delta(\varepsilon)>0$ such that for any sequence $\{[\alpha_{n},\beta_{n}]\}_{n=0}^{\infty}$ of pairwise disjoint subintervals of $I$ satisfying
\bd
\sum\limits_{n=0}^{\infty}|\beta_{n}-\alpha_{n}|<\delta
\ed
it follows
\bd
\sum\limits_{n=0}^{\infty}d(f(\beta_{n}),f(\alpha_{n}))<\varepsilon.
\ed

\edf
\bl\label{l_V}
Let $V$ be defined by (\ref{V}). Then
\bd
\|V\|_{\mathcal{C}^{+}\rightarrow\mathcal{C}^{\pm}}=r_{s}.
\ed
\el
\bpf
The inequality
\bd
\|V\|_{\mathcal{C}^{+}\rightarrow\mathcal{C}^{\pm}}\leq r_{s}
\ed
is trivial. For the equality it is enough to use the constant function $1\in\mathcal{C}^{+}$.
\epf

\bl\label{l_Vy}
Let $V$ be defined by (\ref{V}). Then $Vy\in\mathcal{C}_{A}^{\pm}$ for all $y\in\mathcal{C}^{+}$.
\el
\bpf
For $y\in\mathcal{C}^{+}$ define, according to (\ref{V}), $z\in\mathcal{C}^{\pm}$ as
\bd
z(t):=(Vz)(t)=\left\{
\setlength\arraycolsep{0.1em}\ba{ll}
\ds \int\limits_{0}^{t}y(\sigma)d\sigma,&\;t\in[0,r_{s}]\\
\ds 0,&\;t\in[-\tau,0]
\ea
\right.
\ed
By following Definition \ref{d_AC}, let $\{[\alpha_{n},\beta_{n}]\}_{n=0}^{\infty}$ be any sequence of pairwise disjoint subintervals of $[-\tau,r_{s}]$ satisfying
\bd
\sum\limits_{n=0}^{\infty}|\beta_{n}-\alpha_{n}|<\delta
\ed
for a given $\delta>0$. For a fixed index $n$ it holds
\bd
z(\beta_{n})-z(\alpha_{n})=\left\{
\setlength\arraycolsep{0.1em}\ba{ll}
\ds 0&\textrm{ if }\alpha_{n}<\beta_{n}\leq0\\
\ds \int\limits_{0}^{\beta_{n}}y(\sigma)d\sigma&\textrm{ if }\alpha_{n}\leq0<\beta_{n}\\
\ds \int\limits_{\alpha_{n}}^{\beta_{n}}y(\sigma)d\sigma&\textrm{ if }0<\alpha_{n}<\beta_{n}.\\
\ea
\right.
\ed
Then it easily follows
\bd
|z(\beta_{n})-z(\alpha_{n})|\leq|\beta_{n}-\alpha_{n}|\|y\|_{\mathcal{C}^{+}}
\ed
independently of whether $[\alpha_{n},\beta_{n}]$ falls into $[-\tau,r_{s}]$. Consequently,
\be\label{boundVyAC}
\sum\limits_{n=0}^{\infty}|z(\beta_{n})-z(\alpha_{n})|\leq\sum\limits_{n=0}^{\infty}|\beta_{n}-\alpha_{n}|\|y\|_{\mathcal{C}^{+}}<\delta\|y\|_{\mathcal{C}^{+}}
\ee
and hence $y$ satisfies Definition \ref{d_AC} by setting
\bd
\delta=\frac{\varepsilon}{\|y\|_{\mathcal{C}^{+}}}
\ed
for any $\varepsilon>0$.
\epf

\bl\label{l_GVyAC}
Let $G_{s}$ and $V$ be defined by (\ref{G_{s}}) and (\ref{V}), respectively. Then, under Assumption \ref{coef2}, $G_{s}Vy\in\mathcal{C}_{A}^{+}$ for all $y\in\mathcal{C}^{+}$.
\el
\bpf
For $y\in\mathcal{C}^{+}$ define, according to (\ref{G_{s}}), $z\in\mathcal{C}^{+}$ as
\bd
\setlength\arraycolsep{0.1em}\ba{rcl}
\ds z(t):=(G_{s}Vy)(t)&=&a_{s}(t)(Vy)(t)\\
&&\ds+b_{s}(t)(Vy)(t-\tau)+\int\limits_{-\tau}^{0}c_{s}(t,\theta)(Vy)(t+\theta)d\theta,\;t\in[0,r_{s}].
\ea
\ed
Then, according to (\ref{V}), we have
\bd
z(t)=a_{s}(t)(Vy)(t)+\int\limits_{-t}^{0}c_{s}(t,\theta)(Vy)(t+\theta)d\theta,\;t\in[0,r_{s}],
\ed
whenever $r_{s}<\tau$ while
\bd
z(t)=a_{s}(t)(Vy)(t)+\left\{
\setlength\arraycolsep{0.1em}\ba{ll}
\ds b_{s}(t)(Vy)(t-\tau)+\int\limits_{-\tau}^{0}c_{s}(t,\theta)(Vy)(t+\theta)d\theta,&\;t\in[\tau,r_{s}]\\
\ds\int\limits_{-t}^{0}c_{s}(t,\theta)(Vy)(t+\theta)d\theta,&\;t\in[0,\tau]
\ea
\right.
\ed
whenever $r_{s}\geq\tau$. By following Definition \ref{d_AC}, let $\{[\alpha_{n},\beta_{n}]\}_{n=0}^{\infty}$ be any sequence of pairwise disjoint subintervals of $[0,r_{s}]$ satisfying
\bd
\sum\limits_{n=0}^{\infty}|\beta_{n}-\alpha_{n}|<\delta
\ed
for a given $\delta>0$. Clearly,
\bd
\sum\limits_{n=0}^{\infty}|(\beta_{n}+\theta)-(\alpha_{n}+\theta)|<\delta
\ed
holds for all $\theta\in[-\tau,0]$. It is not difficult to see that, for a fixed index $n$, the bound
\bd
\setlength\arraycolsep{0.1em}\ba{rcl}
\ds|z(\beta_{n})-z(\alpha_{n})|&\leq&\ds|a_{s}(\beta_{n})|\left|(Vy)(\beta_{n})-(Vy)(\alpha_{n})\right|\\
&&\ds+|a_{s}(\beta_{n})-a_{s}(\alpha_{n})|\left|(Vy)(\alpha_{n})\right|\\
&&\ds+|b_{s}(\beta_{n})|\left|(Vy)(\beta_{n}-\tau)-(Vy)(\alpha_{n}-\tau)\right|\\
&&\ds+|b_{s}(\beta_{n})-b_{s}(\alpha_{n})|\left|(Vy)(\alpha_{n}-\tau)\right|\\
&&\ds+\int\limits_{-\tau}^{0}|c_{s}(\beta_{n},\theta)|\left|(Vy)(\beta_{n}+\theta)-(Vy)(\alpha_{n}+\theta)\right|d\theta\\
&&\ds+\int\limits_{-\tau}^{0}|c_{s}(\beta_{n},\theta)-c_{s}(\alpha_{n},\theta)|\left|(Vy)(\alpha_{n}+\theta)\right|d\theta
\ea
\ed
holds thanks to (\ref{V}) both for $r_{s}<\tau$ and $r_{s}\geq\tau$. Now, by applying Lemma \ref{V}, (\ref{boundVyAC}) in Lemma \ref{l_Vy} and by using Assumption \ref{coef2} we get
\bd
\sum\limits_{n=0}^{\infty}|z(\beta_{n})-z(\alpha_{n})|<\delta k\|y\|_{\mathcal{C}^{+}}
\ed
with
\be\label{k}
k:=\big[(\|a_{s}\|_{\mathcal{C}^{+}}+\|b_{s}\|_{\mathcal{C}^{+}}+\|c_{s}\|_{\mathcal{C}^{+},L})+r_{s}(K_{a}+K_{b}+\|K_{c}\|_{L})\big].
\ee
Hence $z$ satisfies Definition \ref{d_AC} by setting
\bd
\delta=\frac{\varepsilon}{k\|y\|_{\mathcal{C}^{+}}}
\ed
for any $\varepsilon>0$.
\epf

\bl\label{l_GVyLIP}
Let $G_{s}$, $V$ and $k$ be defined by (\ref{G_{s}}), (\ref{V}) and (\ref{k}), respectively. Then, under Assumption \ref{coef2}, $G_{s}Vy\in{\rm Lip}_{K}^{+}$ for all $y\in\mathcal{C}^{+}$ with $K\leq k\|y\|_{\mathcal{C}^{+}}$.
\el
\bpf
For $y\in\mathcal{C}^{+}$, Lemma \ref{l_GVyAC} ensures that $G_{s}Vy$ is absolutely continuos in $[0,r_{s}]$, hence also $a.e.$ differentiable. Then it is sufficient to show that $G_{s}Vy$ has a bounded first derivative whenever defined. To this aim observe that for $z=G_{s}Vy$ defined as in the proof of Lemma \ref{l_GVyAC}
\bd
\setlength\arraycolsep{0.1em}\ba{rcl}
\ds z'(t)&=&a_{s}'(t)(Vy)(t)+a_{s}(t)(Vy)'(t)+c_{s}(t,-t)(Vy)(0)\\
&&\ds+\int\limits_{-t}^{0}\left[\frac{\partial c_{s}}{\partial t}(t,\theta)(Vy)(t+\theta)+c_{s}(t,\theta)(Vy)'(t+\theta)\right]d\theta
\ea
\ed
holds for $a.a.$ $t\in[0,r_{s}]$ whenever $r_{s}<\tau$. In fact, under Assumption \ref{coef2}, all coefficients are Lipschitz continuous, hence $a.e.$ differentiable with bounded derivative. Similarly,
\bd
\setlength\arraycolsep{0.1em}\ba{rcl}
\ds z'(t)&=&a_{s}'(t)(Vy)(t)+a_{s}(t)(Vy)'(t)+\left\{
\setlength\arraycolsep{0.1em}\ba{ll}
\ds b_{s}'(t)(Vy)(t-\tau)+b_{s}(t)(Vy)'(t-\tau)\\
\ds c_{s}(t,-t)(Vy)(0)
\ea
\right.\\
&&\ds+\left\{
\setlength\arraycolsep{0.1em}\ba{ll}
\ds\int\limits_{-\tau}^{0}\left[\frac{\partial c_{s}}{\partial t}(t,\theta)(Vy)(t+\theta)+c_{s}(t,\theta)(Vy)'(t+\theta)\right]d\theta,&\;\textrm{ in }[\tau,r_{s}]\\
\ds\int\limits_{-t}^{0}\left[\frac{\partial c_{s}}{\partial t}(t,\theta)(Vy)(t+\theta)+c_{s}(t,\theta)(Vy)'(t+\theta)\right]d\theta,&\;\textrm{ in }[0,\tau]
\ea
\right.
\ea
\ed
holds for $a.a.$ $t\in[0,r_{s}]$ whenever $r_{s}\geq\tau$ and the derivatives of the coefficients are bounded again. Since
\bd
(Vy)'(t)=y(t)
\ed
for all $t\in[0,r_{s}]$ and
\bd
(Vy)(0)=0,
\ed
it follows
\bd
\setlength\arraycolsep{0.1em}\ba{rcl}
\ds z'(t)&=&a_{s}'(t)(Vy)(t)+a_{s}(t)y(t)\\
&&\ds+\int\limits_{-\tau}^{0}\left[\frac{\partial c_{s}}{\partial t}(t,\theta)(Vy)(t+\theta)+c_{s}(t,\theta)y(t+\theta)\right]d\theta
\ea
\ed
for $a.a.$ $t\in[0,r_{s}]$ whenever $r_{s}<\tau$ while
\bd
\setlength\arraycolsep{0.1em}\ba{rcl}
\ds z'(t)&=&a_{s}'(t)(Vy)(t)+a_{s}(t)y(t)\\
&&\ds+\int\limits_{-\tau}^{0}\left[\frac{\partial c_{s}}{\partial t}(t,\theta)(Vy)(t+\theta)+c_{s}(t,\theta)y(t+\theta)\right]d\theta\\
&&\ds+\left\{
\setlength\arraycolsep{0.1em}\ba{ll}
\ds b_{s}'(t)(Vy)(t-\tau)+b_{s}(t)y(t-\tau),&\;\textrm{ in }[\tau,r_{s}]\\
\ds 0,&\;\textrm{ in }[0,\tau]
\ea
\right.
\ea
\ed
for $a.a.$ $t\in[0,r_{s}]$ whenever $r_{s}\geq\tau$. This shows that $z'$ is well-defined in all $[0,r_{s}]$ except for the set of measure zero where the derivative of the coefficients is not defined and also for $t=\tau$ when $r_{s}\geq\tau$. In fact
\bd
z'(\tau^{+})-z'(\tau^{-})=b_{s}(\tau)y(0)\ed
and, in general, $y(0)\neq0$. Nevertheless, for all the other values we have
\bd
|z'(t)|\leq k\|y\|_{\mathcal{C}^{+}}
\ed
with $k$ in (\ref{k}). Since the same bound is valid also for $r_{s}<\tau$, the proof is complete.
\epf

\bl\label{l_I-LG_{s}V}
Let $G_{s}$, $V$ and $\mathcal{L}_{N}^{+}$ be defined by (\ref{G_{s}}), (\ref{V}) and (\ref{LN+}), respectively. Then, under Assumption \ref{coef2},
\be\label{I-LG_{s}Vbound}
\|(I-\mathcal{L}_{N}^{+})G_{s}V\|_{\mathcal{C}^{+}}\leq C\frac{\Lambda_{N}^{+}}{N}
\ee
where $C$ is a constant independent of $N$ and $\Lambda_{N}^{+}$ is the Lebesgue constant relevant to the nodes $\Omega_{N}^{+}$ in (\ref{grid+}). If, in addition, Assumption \ref{cheb+} holds, then
\be\label{I-LG_{s}Vasym}
\lim\limits_{N\rightarrow\infty}\|(I-\mathcal{L}_{N}^{+})G_{s}V\|_{\mathcal{C}^{+}}=0.
\ee
\el
\bpf
By standard interpolation results such as Jackson's Theorem \cite[Theorem 13.3.7]{davis75}, for any given $y\in\mathcal{C}^{+}$, we have
\bd
\setlength\arraycolsep{0.1em}\ba{rcl}
\ds \|(I-\mathcal{L}_{N}^{+})G_{s}Vy\|_{\mathcal{C}^{+}}&\leq&\ds(1+\Lambda_{N}^{+})E_{N}(G_{s}Vy)\\
&&\ds\leq(1+\Lambda_{N}^{+})\left(1+\frac{\pi^{2}}{2}\right)\omega\left(G_{s}Vy;\frac{r_{s}}{2N}\right)\\
&&\ds\leq(1+\Lambda_{N}^{+})\left(1+\frac{\pi^{2}}{2}\right)\frac{r_{s}}{2N}k\|y\|_{\mathcal{C}^{+}}
\ea
\ed
where $E_{N}$ is the best uniform approximation error of $G_{s}Vy$ in $\Pi_{N}^{+}$ and $\omega$, the modulus of continuity, is bounded by virtue of \cite[Theorem 1.5.1]{davis75} and by applying Lemma \ref{l_GVyLIP}. If, in addition, Assumption \ref{cheb+} holds, then (\ref{I-LG_{s}Vasym}) follows by observing that
\bd
\Lambda_{N}^{+}=O(\log{N})
\ed
by virtue of Natanson's Theorem \cite{nat65}.
\epf

\bl\label{l_I-G_{s}Vinv}
Let $G_{s}$ and $V$ be defined by (\ref{G_{s}}) and (\ref{V}), respectively. Then, under Assumption \ref{coef2}, $(I-G_{s}V)^{-1}\in\mathcal{B}(\mathcal{C}^{+})$.
\el
\bpf
We have to prove that for any given $f\in\mathcal{C}^{+}$ there exists a unique $y\in\mathcal{C}^{+}$ solution of
\bd
(I-G_{s}V)y=f.
\ed
By contradiction suppose that, given $f$, there exist $x\neq y$ such that
\bd
(I-G_{s}V)y=f=(I-G_{s}V)x
\ed
and set $z:=y-x$.

If $r_{s}<\tau$, then $z$ solves the Volterra Integral Equation (VIE)
\be\label{VIE1}
z(t)-\int\limits_{0}^{t}k_{s}(t,\sigma)z(\sigma)d\sigma=0,\;t\in[0,r_{s}],
\ee
with
\bd
k_{s}(t,\sigma)=a_{s}(t)+\int\limits_{-t+\sigma}^{0}c_{s}(t,\theta)d\theta,
\ed
as it can be verified by a standard change of integration order for the distributed delay term. The VIE has a continuous kernel $k_{s}$, hence \cite[Theorem 3.12]{kress89} ensures $z=0$.

If $r_{s}\geq\tau$, then $z$ solves the VIE
\be\label{VIE2}
z(t)-\int\limits_{0}^{t}k_{s}(t,\sigma)z(\sigma)d\sigma-\left\{
\setlength\arraycolsep{0.1em}\ba{ll}
\ds \int\limits_{0}^{t-\tau}\left(b_{s}(t)+\int\limits_{-t+\sigma}^{-\tau}c_{s}(t,\theta)d\theta\right)z(\sigma)d\sigma,&\;t\in[\tau,r_{s}]\\
\ds0,&\;t\in[0,\tau].
\ea
\right.
\ee
This latter does not have a continuous kernel due to the presence of the discrete delay term, hence \cite[Theorem 3.12]{kress89} cannot be applied. Nevertheless we can adopt a step-by-step procedure. In fact, being $r_{s}\geq\tau$, there exist an integer $\alpha$ and a $\beta\in[0,\tau)$ such that $r_{s}=\alpha\tau+\beta$. We prove that $z(t)=0$ for all $t\in[(i-1)\tau,i\tau]$, $i=1,\ldots,\alpha$, and also for $t\in[\alpha\tau,r]$. For $i=1$ we have $t\in[0,\tau]$ and the VIE (\ref{VIE2}) reduces to (\ref{VIE1}) for which \cite[Theorem 3.12]{kress89} provides $z(t)=0$ for all $t\in[0,\tau]$. For $i=2$ we have $t\in[\tau,2\tau]$ and the last integral in (\ref{VIE2}) over $[0,t-\tau]$ involves the solution over $[0,\tau]$, which has just been proven to be null, showing that $z(t)=0$ for all $t\in[\tau,2\tau]$ as well. Iteration of the process leads to $z(t)=0$ for all $t\in[0,r_{s}]$ since the same argument holds for the last interval $[\alpha\tau,r]$, too.

Eventually, independently of whether $r_{s}<\tau$ or $r_{s}\geq\tau$, we get $z=0$ and $y=x$, which contradicts the initial assumption.
\epf

\bl\label{l_I-LG_{s}Vinv}
Let $G_{s}$, $V$ and $\mathcal{L}_{N}^{+}$ be defined by (\ref{G_{s}}), (\ref{V}) and (\ref{LN+}), respectively. Then, under Assumptions \ref{cheb+} and \ref{coef2} and for sufficiently large $N$, $(I-\mathcal{L}_{N}^{+}G_{s}V)^{-1}\in\mathcal{B}(\mathcal{C}^{+})$. In particular,
\bd
\|(I-\mathcal{L}_{N}^{+}G_{s}V)^{-1}\|_{\mathcal{C}^{+}}\leq2\|(I-G_{s}V)^{-1}\|_{\mathcal{C}^{+}}.
\ed
\el
\bpf
The thesis follows immediately by observing that
\bd
I-\mathcal{L}_{N}^{+}G_{s}V=(I-G_{s}V)+(I-\mathcal{L}_{N}^{+})G_{s}V
\ed
and by applying the Banach's perturbation Lemma \cite[Theorem 10.1]{kress89} whose hypotheses are satisfied thanks to Lemmas \ref{l_I-LG_{s}V} and \ref{l_I-G_{s}Vinv}.
\epf

\begin{thebibliography}{10}

\bibitem{ans71}
P.~M. Anselone.
\newblock {\em Collectively compact operator approximation theory}.
\newblock Prentice-{H}all, Englewood {C}liffs, {N}ew {J}ersey, 1971.

\bibitem{belzen03}
A.~Bellen and M.~Zennaro.
\newblock {\em Numerical methods for delay differential equations}.
\newblock Numerical Mathemathics and Scientifing Computing series. Oxford
  University Press, 2003.

\bibitem{bddm92}
A.~Bensoussan, G.~Da~Prato, M.~C. Delfour, and S.~K. Mitter.
\newblock {\em Representation and control of infinite dimensional systems},
  volume I and II.
\newblock Birkh{\"{a}}user, 1992, 1993.

\bibitem{bt69}
J.~G. Borisovi{\u{c}} and A.~S. Turbabin.
\newblock On the {C}auchy problem for linear nonhomogeneous differential
  equations with retarded arguments.
\newblock {\em Dokl. Akad. Nauk SSSR}, 185(4):741--744, 1969.
\newblock English transl. Soviet Math. Dokl., 10(2):401-405, 1969.

\bibitem{breth04}
D.~Breda.
\newblock {\em Numerical computation of characteristic roots for delay
  differential equations}.
\newblock PhD thesis, PhD in Computational Mathematics, Universit\`a di Padova,
  2004.

\bibitem{brean206}
D.~Breda.
\newblock Solution operator approximation for characteristic roots of delay
  differential equations.
\newblock {\em Appl. Numer. Math.}, 56(3-4):305--317, 2006.

\bibitem{brejmaa09}
D.~Breda.
\newblock Nonautonomous delay differential equations in {H}ilbert spaces and
  {L}yapunov exponents.
\newblock 2009.
\newblock Preprint submitted for publication.

\bibitem{bresisc05}
D.~Breda, S.~Maset, and R.~Vermiglio.
\newblock Pseudospectral differencing methods for characteristic roots of delay
  differential equations.
\newblock {\em SIAM J. Sci. Comput.}, 27(2):482--495, 2005.

\bibitem{breifac06}
D.~Breda, S.~Maset, and R.~Vermiglio.
\newblock Numerical computation of characteristic multipliers for linear time
  periodic delay differential equations.
\newblock In C.~Manes and P.~Pepe, editors, {\em Time Delay Systems 2006},
  volume~6 of {\em IFAC Proceedings Volumes}. Elsevier, 2006.

\bibitem{brenm09}
D.~Breda, S.~Maset, and R.~Vermiglio.
\newblock Numerical approximation of characteristic values of partial retarded
  functional differential equations.
\newblock {\em Numer. Math.}, 113(2):181--242, 2009.

\bibitem{breifac210}
D.~Breda, S.~Maset, and R.~Vermiglio.
\newblock On discretizing the semigroup of solution operators for linear time
  invariant - time delay systems.
\newblock In {\em Time Delay Systems 2010}, IFAC Proceedings Volumes, 2010.
\newblock Submitted for publication.

\bibitem{bmbas04}
E.~A. Butcher, H.~T. Ma, E.~Bueler, V.~Averina, and Z.~Szabo.
\newblock Stability of linear time-periodic delay-differential equations via
  chebyshev polynomials.
\newblock {\em Int. J. Numer. Meth. Engng}, 59:895--922, 2004.

\bibitem{chat83}
F.~Chatelin.
\newblock {\em Spectral approximation of linear operators}.
\newblock Academic Press, New York, 1983.

\bibitem{davis75}
P.~J. Davis.
\newblock {\em Interpolation \& approximation}.
\newblock Dover, New York, USA, 1975.

\bibitem{del77}
M.~C. Delfour.
\newblock State theory of linear hereditary differential systems.
\newblock {\em J. Differ. Equations}, 60:8--35, 1977.

\bibitem{dm72}
M.~C. Delfour and S.~K. Mitter.
\newblock Hereditary differential systems with constant delays. {I}. {G}eneral
  case.
\newblock {\em J. Differ. Equations}, 12:213--235, 1972.

\bibitem{bbck07}
M.~di~Bernardo, C.~J. Budd, A.~R. Champneys, and P.~Kowalczyk.
\newblock {\em Piecewise-smooth dynamical systems: theory and applications}.
\newblock Number 163 in AMS series. Springer Verlag, New York, USA, 2007.

\bibitem{diekmann95}
O.~Diekmann, S.~A. van Gils, S.~M. Verduyn~Lunel, and H.~O. Walther.
\newblock {\em Delay Equations - Functional, Complex and Nonlinear Analysis}.
\newblock Number 110 in AMS series. Springer Verlag, New York, USA, 1995.

\bibitem{engacm02}
K.~Engelborghs, T.~Luzyanina, and D.~Roose.
\newblock Numerical bifurcation analysis of delay differential equations using
  {DDE-BIFTOOL}.
\newblock {\em ACM T. Math. Software}, 28(1):1--21, 2002.

\bibitem{engsjna02}
K.~Engelborghs and D.~Roose.
\newblock On stability of {LMS} methods and characteristic roots of delay
  differential equations.
\newblock {\em SIAM J. Numer. Anal.}, 40(2):629--650, 2002.

\bibitem{fab14}
G.~Faber.
\newblock {\"{U}}ber die interpolatorische darstellung stetiger funktionen.
\newblock {\em Jahresber. Deut. Math. Verein.}, 23:192--210, 1914.

\bibitem{farmer82}
D.~Farmer.
\newblock Chaotic attractors of an infinite-dimensional dynamical system.
\newblock {\em Physica D}, 4:605--617, 1982.

\bibitem{gp06}
D.~E. Gilsinn and F.~A. Potra.
\newblock Integral operators and delay differential equations.
\newblock {\em J. Integral Equations Appl.}, 94(6):297--336, 2006.

\bibitem{ggk90}
I.~Gohberg, S.~Goldberg, and M.~A. Kaashoek.
\newblock {\em Classes of linear operators}.
\newblock Number~49 in Operator Theory: Advances and Applications. Birkhauser,
  New York, USA, 1990.

\bibitem{hrs08}
S.~Hadd, A.~Rhandi, and R.~Schnaubelt.
\newblock Feedback theory for time-varying regular linear systems with input
  and state delays.
\newblock {\em IMA J. Math. Control Inform.}, 25(1):85--110, 2008.

\bibitem{hvl93}
J.~K. Hale and S.~M. Verduyn~Lunel.
\newblock {\em Introduction to functional differential equations}.
\newblock Number~99 in AMS series. Springer Verlag, New York, USA, 2nd edition,
  1993.

\bibitem{ig02}
T.~Insperger and G.~St\'ep\'an.
\newblock Semi-discretization method for delayed systems.
\newblock {\em Int. J. Numer. Meth. Engng}, 55:503--518, 2002.

\bibitem{ik91}
K.~Ito and F.~Kappel.
\newblock A uniformly differentiable approximation scheme for delay systems
  using splines.
\newblock {\em Appl. Math. Opt.}, 23:217--262, 1991.

\bibitem{jar08}
E.~Jarlebring.
\newblock {\em The spectrum of delay-differential equations: numerical methods,
  stability and perturbation}.
\newblock PhD thesis, Inst. Comp. Math, TU Braunschweig, 2008.

\bibitem{kap86}
F.~Kappel.
\newblock {\em Semigroups and delay equations}.
\newblock Number 152 ({T}rieste, 1984) in Pitman Res. Notes Math. Ser. Longman
  Sci. Tech., Harlow, 1986.

\bibitem{kra59}
N.~Krasovskii.
\newblock {\em Stability of Motion}.
\newblock Moscow, 1959.
\newblock English transl. Stanford University Press, 1963.

\bibitem{kress89}
R.~Kress.
\newblock {\em Linear integral equations}.
\newblock Number~82 in AMS series. Springer-Verlag, New York, USA, 1989.

\bibitem{kry56}
V.~I. Krylov.
\newblock Convergence of algebraic interpolation with respect to roots of
  {C}hebyshev's polynomial for absolutely continuous functions of bounded
  variation.
\newblock {\em Dokl. Akad. Nauk SSSR}, 107:362--365, 1956.

\bibitem{nat65}
I.~P. Natanson.
\newblock {\em Constructive function theory Vol. {III}}.
\newblock Frederick Ungar Publ., New York, USA, 1965.

\bibitem{pei82}
G.~H. Peichl.
\newblock A kind of ``history space" for retarded functional differential
  equations and representation of solutions.
\newblock {\em Funkc. Ekvacioj-SER I}, 25:245--256, 1982.

\bibitem{tref00}
L.~N. Trefethen.
\newblock {\em Spectral methods in {MATLAB}}.
\newblock Software - Environment - Tools series. SIAM, Philadelphia, USA, 2000.

\bibitem{vlr08}
K.~Verheyden, T.~Luzyanina, and D.~Roose.
\newblock Efficient computation of characteristic roots of delay differential
  equations using lms methods.
\newblock {\em J. Comput. Appl. Math.}, 214(1):209--226, 2008.

\bibitem{vin78}
R.~B. Vinter.
\newblock On the evolution of the state of linear differential delay equations
  in {$M^{2}$}: properties of the generator.
\newblock {\em J. Inst. Maths. Applics.}, 21:13--23, 1978.

\bibitem{vz09}
T.~Vyhl{\'\i}dal and P.~Z{\'\i}tek.
\newblock Mapping based algorithm for large-scale computation of
  quasi-polynomial zeros.
\newblock {\em IEEE T. Automat. Cont.}, 54(1):171--177, 2009.

\bibitem{wu96}
J.~Wu.
\newblock {\em Theory and applications of partial functional differential
  equations}.
\newblock Number 119 in AMS series. Springer-Verlag, New York, USA, 1996.

\end{thebibliography}

\end{document}